\documentclass[11pt]{amsart}

\usepackage{url}
\usepackage{psfrag}
\usepackage{graphicx}
\usepackage{comment}
\usepackage{array}
\usepackage{amsmath,amsfonts,amssymb,amsxtra,amsthm}
\usepackage{mathrsfs}
\usepackage{lineno}
\usepackage{verbatim}
\usepackage{xy}
\usepackage{times}
\xyoption{all}
\usepackage{stmaryrd}
\usepackage{appendix}
\usepackage{xr}

\DeclareMathOperator{\betti}{b_1}

\DeclareMathOperator{\Cyl}{Cyl}

\DeclareMathOperator{\edgegraphs}{EdgeGraphs}

\DeclareMathOperator{\img}{Im}

\DeclareMathOperator{\link}{lk}

\DeclareMathOperator{\rk}{rk}

\DeclareMathOperator{\scott}{sc}

\DeclareMathOperator{\st}{St}

\DeclareMathOperator{\vertexgraphs}{VertexGraphs}

\def\cw{\textsc{CW}}

\newcommand{\E}{\mathcal{E}}

\newcommand{\RR}{\mathcal{R}}

\newcommand{\YY}{\mathcal{Y}}
\newcommand{\ZZ}{\mathcal{Z}}

\def\xtilde{\ensuremath{\widetilde{{X}}}}

\def\zbar{\ensuremath{\overline{Z}}}

\def\sl2c{\ensuremath{{SL}(2,\mathbb{C})}}
\def\t1sl2c{{\mathfrak sl}_2\mathbb{C}}
\def\free{\ensuremath{\mathbb{F}}}
\def\zee{\mathbb{Z}}

\def\immerses{\looparrowright}
\def\onto{\twoheadrightarrow}
\def\into{\hookrightarrow}

\newcommand{\adjoin}[2]{\ensuremath{#1\!\left[#2\right]}}
\newcommand{\adjoinroot}[3]{\ensuremath{\adjoin{#1}{\sqrt[#3]{#2}}}}
\newcommand{\grushko}[3]{\ensuremath{#1_1*\cdots*#1_{#2}*\free_{#3}}}

\makeatletter
\newcommand{\Doubletwo}[3]{ \left\{ #1 \mid #3\right\} }
\newcommand{\Doubleone}[1]{ \left\{ #1 \right\} }
\newcommand{\Doublethr}[3]{ \left\{ #1 \right\}_{#3} }
\newcommand{\set}[1]{%
\@ifnextchar:{\Doubletwo{#1}}{\@ifnextchar_{\Doublethr{#1}}{\Doubleone{#1}}}%
}
\makeatother

\makeatletter
\newcommand{\grouptwo}[3]{ \langle #1 \mid #3\rangle }
\newcommand{\groupone}[1]{ \langle #1 \rangle }
\newcommand{\group}[1]{%
\@ifnextchar:{\grouptwo{#1}}{\groupone{#1}}%
}
\makeatother

\def\sonecomp{\Gamma_{\!\circ}}
\def\underlying{\ensuremath{\Gamma_{\!\mathrm{U}}}}
\def\compcomp{\Gamma_{\!\infty}}

\newcommand{\term}[1]{{\emph{#1}}}
\newcommand{\scare}[1]{`#1'}

\newcommand{\mobius}{M\"obius}

\newtheorem{theorem}{Theorem}[section]
\newtheorem{lemma}[theorem]{Lemma}
\newtheorem{example}[theorem]{Example}
\newtheorem{corollary}[theorem]{Corollary}

\theoremstyle{definition}
\newtheorem{definition}[theorem]{Definition}

\theoremstyle{remark}
\newtheorem{remark}[theorem]{Remark}
\newtheorem{convention}[theorem]{Convention}

\newcommand{\mnote}[1]{}

\title[Scott complexity and adjoining roots]{Krull dimension for limit
  groups {III}:\\
  Scott complexity and adjoining roots to finitely generated groups}

\author{Larsen Louder}
\address{Department of Mathematics\\
University of Michigan \\
Ann Arbor, MI 48109-1043\\
USA}
\email[Larsen Louder]{llouder@umich.edu, lars@d503.net}

\keywords{limit groups, krull dimension, Stallings's folding}

\subjclass[2000]{Primary: 20F65; Secondary: 20E05, 20E06}

\thanks{Most of this research was done while at the University of
  Utah. The author also gratefully acknowledges support from the
  National Science Foundation.}

\date{\today}

\begin{document}

\begin{abstract}
  This is the third paper in a sequence on Krull dimension for limit
  groups, answering a question of Z.~Sela. We give generalizations of
  the well known fact that a nontrivial commutator in a free group is
  not a proper power to both graphs of free groups over cyclic
  subgroups and freely decomposable groups.  None of the paper is
  specifically about limit groups.
\end{abstract}

\maketitle


\section{Introduction}

\thispagestyle{empty}

\par This and the companion paper~\cite{louder::roots} contain an
analysis of sequences of limit groups obtained through a certain
systematic process of adjoining roots and passing to limit group
quotients. An analysis of arbitrary finitely generated groups obtained
through adjoining roots is, of course, not practical, however
something can be said if only limit groups are considered.

\par This paper is preceded by~\cite{louder::strict,louder::stable},
and is followed up by~\cite{louder::roots}. In~\cite{louder::stable}
we reduce the problem of finite Krull dimension to a statement about
adjoining roots to limit groups. The argument needed divides neatly
into two parts. The first, dealt with in this paper, is a study of
groups obtained by adjoining roots to almost all finitely generated
groups. We handle the case of arbitrary limit groups in the last.

\begin{definition}[\cite{scottcoherent,swarup::scott}]
  Let $G$ be a finitely generated group with Grushko decomposition
  $\grushko{G}{p}{q}$. The \term{Scott complexity} of $G$ is the
  ordered pair $\scott(G)=(q-1,p)$.
\end{definition}

Let $G$ be a group. An element $g\in G$ is \term{indivisible} if isn't
a proper power, i.e., if $h^k=g$ then $k=\pm 1$. The group $\langle
G,\gamma'_i\vert(\gamma'_i)^{k_i} = \gamma_i,i=1\dotsc n\rangle$ is
denoted $\adjoinroot{G}{\gamma_i}{k_i}$. For notational convenience we
usually suppress the ``$k_i$'' from the notation. This is justified by
the fact that the degree of the root is incidental as long as
$\gamma_i$ indivisible in $G,$ as we see in
Lemma~\ref{addroottofreegroup}. The letter $\free$ represents a
nonabelian free group. The \term{corank} of a finitely generated group
is the maximal rank of a free group it maps onto.

\begin{theorem}[Scott complexity and addjoining roots to groups]
  \label{maintheorem}
  Suppose that $\phi\colon G\into H$ and $H$ is a quotient of
  $G'=\adjoinroot{G}{\gamma_i}{k_i},$ $\boldsymbol{\gamma}_i$ a
  collection of distinct conjugacy classes of indivisible elements of
  $G$ such that $\boldsymbol{\gamma}_i\neq\boldsymbol{\gamma}_j^{-1}$
  for all $i,j$ and $\gamma_i\in\boldsymbol{\gamma}_i$. Then
  $\scott(G)\geq\scott(H)$. If equality holds and $H$ has no $\zee_2$
  free factors, there are presentations of $G$ and $H$ as
  \[
  G\cong\grushko{G}{p}{q}^G,\quad H\cong\grushko{H}{p}{q}^H
  \]
  a partition of $\set{\boldsymbol{\gamma}_i}$ into subsets
  $\boldsymbol{\gamma}_{j,i},$ $j=0,\dotsc,p,$ $i=1,\dotsc,i_p,$
  representatives $\gamma_{j,i}\in G_j,\boldsymbol{\gamma}_{j,i},$
  $i\geq 1,$ $\gamma_{0,i}\in\free_q^G,\boldsymbol{\gamma}_{0,i},$
  such that with respect to the presentations of $G$ and $H$:
  \begin{itemize}
    \item $\phi(G_i)<H_i$
    \item $\adjoinroot{G_j}{\gamma_{j,i}}{k_{j,i}}\onto H_j$
    \item $\phi(\free_q^G)<\free_q^H$
    \item
    $\free^G_q=\langle\gamma_{0,1}\rangle*\dotsb*\langle\gamma_{0,i_0}\rangle*F$
    \item $\free^H_q=\langle\sqrt{\gamma_{0,1}}\rangle*\dotsb*\langle\sqrt{\gamma_{0,i_0}}\rangle*F$
    \item
    $G'\cong\adjoinroot{G_1}{\gamma_{1,i}}{}*\dotsb*\adjoinroot{G_p}{\gamma_{p,i}}{}*\langle\sqrt{\gamma_{0,1}}\rangle*\dotsb*\langle\sqrt{\gamma_{0,i_0}}\rangle*F$
  \end{itemize}
  All homomorphisms are those suggested by the presentations.
\end{theorem}

\begin{example}[No $\zee_2$ free factors of $H$ is a necessary hypothesis]
  Let $G=G_1*G_2,$ $H=H_1*\zee_2,$ where
  $H_1=G_1*G_2/\langle\langle\alpha_1=\alpha_2\rangle\rangle,$
  $\alpha_i\neq 1\in G_i,$ and $x$ generates $\zee_2$. Then $H$ is a
  quotient of $\adjoinroot{G}{\alpha_1\alpha_2}{}$. The inclusion
  $G\into H$ maps $G_1$ to $G_1<H_1$ and $G_2$ to $xG_2x<xH_1x$. Then
  $\alpha_1\to\alpha_1$ and $\alpha_2\to x\alpha_1x$ and
  $\alpha_1\alpha_2\to(\alpha_1x)^2$.
\end{example}

\par In the case $\scott(G)=\scott(H)=(q-1,0),$ and adjoining a root
to a single nontrivial element,~Theorem~\ref{maintheorem} reduces to a
theorem of Baumslag~\cite{baumslagpowers}. The starting point for our
approach is a theorem of Shenitzer~\cite{shenitzer} which states that
if an amalgamation $\free_n*_{\langle t\rangle}\free_m$ is free
(necessarily of rank $n+m-1$!), then $t$ is a basis element in at
least one vertex group. If $\adjoinroot{G}{\gamma}{k}$ is free then
Theorem~\ref{maintheorem} reduces to Shenitzer's theorem: $\gamma$ is
not a basis element of $\langle\sqrt[k]{\gamma}\rangle,$ therefore
$\gamma$ is a basis element of $F$.

\par A modern proof of Shenitzer's Theorem goes as follows: If
$\free_n*_{\langle t\rangle}\free_m$ is free, then each map of a
vertex group into $\free_n*_{\langle t\rangle}\free_m$ may be
represented by an immersions $\iota_{n\vert m}\colon\Gamma_{n\vert
  m}\immerses R_{n+m-1}$. Represent $t$ by immersions of
$S^1\immerses\Gamma_{n\vert m},$ and build a graph of spaces $X$ by
gluing the boundary components of an annulus to
$\Gamma_n\sqcup\Gamma_m$ along the immersions of $S^1$. Extend
$\iota_{n\vert m}$ to a map $\eta\colon X\to R_{n+m-1}$. Pull back
midpoints of edges of $\mathrm{R}_{n+m-1}$ to produce embedded graphs
in $X$ transverse to both $\Gamma_{n}$ and $\Gamma_m$.  The preimage
graphs must be forests, otherwise $\free_n*_{\langle t
  \rangle}\free_m\onto\free_{n+m-1}$ has nontrivial kernel. This
implies that the representation of $t$ as an immersion, in at least
one of $\Gamma_{n\vert m},$ must cover some edge only one time, i.e.,
it is a basis element in one of the factors.

\par If we drop the hypothesis that $\free_n*_{\langle t
  \rangle}\free_m$ is free the means to conclude that preimages of
midpoints of edges are forests disappears. The next theorem shows that
this hypothesis can be weakened.

\begin{theorem}[Improved Shenitzer's Theorem]
  \label{improvedshenitzer}
  Let $G=F_1*_{\langle t\rangle}F_2,$ $F_1$ and $F_2$ free of rank $n$
  and $m,$ respectively. If $t$ is indivisible in at least one vertex
  group and $G$ has corank ${n+m-1},$ then $G$ is free.
\end{theorem}

\par Our approach to Theorem~\ref{maintheorem} is to carry out a more
careful analysis of the space constructed to prove Shenitzer's
theorem. Before moving on to this, we give a short outline of some
existing work on ranks of subgroups generated by solutions to
equations defined over $\free$.

\par The first equation over the free group to receive much attention
was Vaught's equation $\Omega=\set{a^2b^2c^2=1}$. Lyndon showed that if
$a^2b^2c^2=1$ in a free group, then $a,$ $b,$ and $c$ commute. This
characterization of solutions amounts to the fact that a commutator in
a free group isn't a square. Lyndon and
Sch{\"u}tzenberger~\cite{lyndonabc} extended this fact to equations
$\Omega=\set{a^pb^qc^r=1},$ $p,q,r\geq 2,$ of the same type.

\par More generally, Baumslag showed that if $(y_1,\dotsc,y_{n+1})$ is
a solution to $\omega(x_1,\ldots,x_n)=x_{n+1}^k,$ $k>1$ in $\free,$
and $\omega$ is neither a proper power nor a basis element in $\langle
x_1,\dotsc,x_n\rangle,$ then the subgroup of $\free$ generated by
$\set{y_1,\dotsc,y_{n+1}}$ has rank less than $n$. 

\par In a remark, Baumslag seems to suggest a conjugacy separability
problem for elements of the free group: If $\alpha$ and $\beta$ are
nonconjugate elements of $\free_n,$ then either (without loss) there
is a free factorization $\free_n\cong F*\langle\alpha'\rangle$ such
that $\alpha$ is a power of $\alpha'$ and $\beta$ is conjugate into
$F,$ or $\langle\free_n,t\vert t\alpha t^{-1}\beta^{-1}\rangle$ has
corank less than $n$. The next corollary is a resolution of this
question.

\begin{corollary}
  \label{maximumrankconjugacy} 
  Let $F$ be a free group of rank $n>1,$ and let $\mathrm{Z}_i<F$ be
  finitely many distinct conjugacy classes of maximal cyclic
  subgroups of $F$. Let $t_j$ be stable letters,
  $\gamma^1_j,\gamma^2_j$ elements of
  $\bigcup_i\mathrm{Z}_i\setminus\set{1}$. Let $G$ be the group
  $F*\langle t_j\rangle/\langle\!\langle t_j\gamma^1_j
  t_j^{-1}=\gamma^2_j\rangle\!\rangle$. Let $\sim$ be the equivalence
  relation generated by
  \[\mbox{$\mathrm{Z}_i\sim\mathrm{Z}_{i'}$ $\iff$ There exists $j$ such that  $\gamma^1_j\in\mathrm{Z}_i$ and $\gamma^2_j\in\mathrm{Z}_{i'}$}\]

  If $\sim$ has no singleton equivalence classes and $G$ has corank $n,$
  then for some $i$ there is a free factorization of $F$ as $F\cong
  F_1*\mathrm{Z}_i$ such that every element $\gamma^s_j$ is conjugate
  into either $F_1$ or $\mathrm{Z}_i$ and $F_1$ contains a conjugate
  of some $\gamma^s_j$.
\end{corollary}

This corollary should be compared to the criterion for malnormality of
rank two subgroups of free groups given in~\cite{fmr02}. Say that a
subgroup $H$ of a free group $\free$ is \term{isolated} if it is
closed under taking roots, \term{isolated on generators} if it is
closed under taking roots of generators, and \term{malnormal on
  generators} if, for all $g\in \free\setminus H$ and basis elements
$h$ of $H,$ $ghg^-1\notin H$. Fine, Myasnikov, and Rosenberger prove
that if a rank two subgroup of a free group is isolated and malnormal
on generators, then it is malnormal. By Theorem~\ref{maintheorem}, a
rank two subgroup of a free group is isolated if and only if it is
isolated on generators, and, if not, is not malnormal on
generators. If malnormality on generators fails for some reason other
than isolation, by Corollary~\ref{maximumrankconjugacy}, any two
elements nonconjugate in $H$ but conjugate in $\free$ must be
conjugate to a basis of $H$. A non-malnormal subgroup must contain such
a pair of elements, and so malnormality is implied by malnormality on
generators. The author suspects that the paper of Fine et al.\ contains
a proof of Corollary~\ref{maximumrankconjugacy} for the case
considered in this paragraph.


\subsection*{Acknowledgments}

It's hard to overstate my gratitude to Mark Feighn for his careful
reading of more versions of this paper than I can count. Many thanks
must also go to Mladen Bestvina for pointing out the proof of
Shenitzer's theorem and teaching me how to fold, many moons ago.



\section{Injections, immersions, and graphs of spaces}

\par To analyze homomorphisms of the type
in~\ref{maintheorem},~\ref{improvedshenitzer},
and~\ref{maximumrankconjugacy}, we need to construct spaces which
efficiently represent injections of groups. Given an injection
$F_1\into F_2$ of free groups, Stallings constructs graphs $\Gamma_1$
and $\Gamma_2$ and a map $\Gamma_1\to\Gamma_2$ which, under suitable
identifications of $\pi_1(\Gamma_{1,2})$ with $F_{1,2},$ represents
given homomorphism. In this section we generalize his construction to
spaces which represent injections of groups which aren't necessarily
freely indecomposable, but which are strong enough to promote
Stallings' type results from free groups to Grushko free
factorizations of freely decomposable groups. Our maps give absolutely
no information about restrictions to freely indecomposable free
factors.

\par We give a brief review of immersions before we generalize
immersions to ``relative graphs,'' or relative one-complexes, a class
of spaces slightly larger than that of graphs.

\subsection{Immersions}

\par A map of graphs $\varphi\colon\Gamma_1\to\Gamma_2$ is
\term{combinatorial} if
\begin{itemize}
\item $\varphi(\Gamma_1^{(0)})\subset\Gamma_2^{(0)}$
\item If $e$ is an edge of $\Gamma_1$ then there is an edge $f$ of
  $\Gamma_2$ so that $\varphi\vert_{e^{\circ}}\colon e^{\circ}\to
  f^{\circ}$ is a homeomorphism.
\end{itemize}
A combinatorial map induces, for each vertex $v$ of $\Gamma,$ a map
$\varphi_v\colon\link(v)\to\link(\varphi(v))$. If each $\varphi_v$ is
injective then $\varphi$ is an \term{immersion}. Our goal is to
represent an injection $G=\grushko{G}{p_G}{q_G}\into
H=\grushko{H}{p_H}{q_H}$ as an \scare{immersion} of suitable cell
complexes with fundamental groups $G$ and $H$. We first translate the
link condition into \term{local relative $\pi_1$ injectivity}, then we
extend the translation to relative graphs.

\par Let $\varphi\colon \Gamma_1\to\Gamma_2$ be combinatorial, $\pi_1$
injective, and suppose that the stars of vertices of $\Gamma_2$ are
embedded. For each vertex $v$ there is an induced map
\[\pi_1(\varphi_v)\colon\pi_1(\st(v),\partial\st(v))
\to\pi_1(\st(\varphi(v)),\partial\st(\varphi(v)))\] If $\varphi$ is an
immersion then $\pi_1(\varphi_v)$ is injective for all vertices $v$ of
all combinatorial representatives of the topological realization of
$\varphi$.


\par If $\Box(X)$ is a collection of subspaces of a space $X$ denote
the disjoint union of elements of $\Box(X)$ by $\Box_{\bullet}(X)$.

\begin{definition}[Relative Graph]
  \par A \term{relative graph} is a topological space $W$ with the
  structure of a \cw-pair of cell complexes $(X,\YY(X))$ where
  \begin{itemize}
  \item The topological realization of $X$ is homeomorphic to $W$.
  \item $\YY(X)$ is a collection of connected subcomplexes of
    $X$. If $Y$ and $Y'\in\YY(X)$ meet, then they meet in
    finitely many vertices. If it's clear which relative graph we're
    referring to, the `$(X)$' of `$\YY(X)$' will be
    suppressed.
  \item If $e$ is a cell not contained in $\YY_{\bullet}$
    then $e$ is at most one-dimensional
  \item There is a (not necessarily connected) graph $\Gamma_X,$ and
    for each valence one or zero vertex $v$ of $\Gamma_X,$ there are
    finitely many maps $\eta_{v,i}\colon v\to \YY_{\bullet}(X)$
  \item $W$ is homeomorphic to the quotient space
    $\Gamma_X\sqcup\YY_{\bullet}(X)/v\sim\eta_{v,i}(v)$.
  \item Each $Y\in\YY(X)$ has nontrivial fundamental group.
  \end{itemize}
\end{definition}

\par Let $(X,\YY(X))$ be a relative graph. The \term{zero
  skeleton} of $(X,\YY(X))$ is the set of zero cells not
contained in any element of $\YY(X),$ along with the connected
components of the union of the elements of $\YY(X)$ in $X$.
Pick an element $v$ of the zero skeleton, let $e_1,\dotsc,e_n$ be the
\term{oriented} edges of $\Gamma_X$ such that $\tau(e_i)\in v,$ and
identify the edges $e_j$ with intervals $I_1,\dotsc,I_n$. Then the
star of $v$ is the complex $\st(v)=v\cup\sqcup I_j/(1\in
I_j\sim\tau(e_j))$. There is a map $\st(v)\to X$ such that $v\to v$ is
the identity map and $I_j\to e_j$ is simply the prescribed
identification.

\par The \term{essential zero skeleton} of a relative graph $X,$ $X^E$
is obtained by removing all valence two vertices from the zero
skeleton.
The \term{inessential zero skeleton} consists of the valence two
vertices of $\Gamma_X,$ and is denoted $X^I$.

A relative graph $X$ is \term{admissible} if every map $\st(v)\to X$
is an embedding and, for all $v,w\in X^E,$ $\st(v)$ and $\st(w)$ don't
share any edges of $\Gamma_X$.

\begin{definition}
  Let $X$ and $X'$ be admissible relative graphs. A map $\varphi\colon
  X\to X'$ is \term{combinatorial} if the preimage of the interior of
  an edge of $X'$ is the union of interiors of edges from $X$. The
  restriction of $\varphi$ to the interior of an edge of $X$ is a
  homeomorphism onto the corresponding interior in $X'$.
\end{definition}

\begin{definition}
  \label{immersion}
  Let $\varphi\colon X'\to X$ be a combinatorial map of relative
  graphs, and suppose that $(X,\YY(X))$ is admissible. We say
  that $\varphi$ is an \term{immersion} if neither of the following
  conditions holds:
\begin{itemize}
  \item If, for some
    vertex $v\in\Gamma_X\setminus\YY_{\bullet}(X)$ and
    $w\in\varphi^{-1}(v),$
    \[\pi_1(\varphi_w)\colon\pi_1(\st(w),\partial\st(w))\to\pi_1(\st(v),\partial\st(v))\]
    is not injective then $\varphi$ is not an immersion.
  \item Suppose $\pi_1(\varphi_w)$ is injective for all such $w$ and
    $v$ as in $1$. If, for some $Y\in\YY(X)$ and connected
    component $N$ of $\varphi^{-1}(\st(Y)),$ the map
    \[\pi_1(\varphi_N)\colon\pi_1(N,\partial
    N)\to\pi_1(\st(Y),\partial\st(Y))\] is not injective then
    $\varphi$ is not an immersion.
\end{itemize}
\end{definition}

\par Definition~\ref{immersion} requires some explanation. Suppose
$\varphi$ doesn't satisfy the first bullet, i.e., it is an immersion
at ordinary vertices. Now consider $Y\in\YY(X),$ $\st(Y),$ and
a connected component $N$ of $\varphi^{-1}(\st(Y))$. The preimage $N$
is the union of elements $Y'_i\in\YY(X')$ and edges of
$\Gamma_{X'}$. If $v\in\varphi^{-1}(\partial\st(Y))\cap Z$ then at
most one edge of $\Gamma_{X'}\cap N$ meets $v,$ since if there were
two, then the map on $\st(v)$ wouldn't be injective, contrary to
hypothesis. Thus $N=\st(V)$ for some subcomplex $V$ of $X',$ and
$\partial(N)=\varphi^{-1}(\partial\st(Y))\cap N$. Thus the second
condition has the same form as the first in the event that the first
doesn't hold. Stallings calls paths representing elements of
$\ker(\pi_1(\varphi_N))$ \emph{binding ties}~\cite{stallingsgrushko}.

An \term{edge path} in a relative graph is a combinatorial map of a
subdivided interval. An edge path which is an immersion is
\term{reduced}. An edge path is homotopic, relative to its endpoints,
to a reduced edge path, and any two reduced edge paths homotopic
relative to endpoints are equivalent, as combinatorial objects, via
homotopies supported on those subsegments of the interval with image
contained in some $\YY_{\bullet}(X)\subset X$.

\mnote{this is screwey now}

\par A \term{fold} of a relative graph $X$ is a map of the following
type: Let $e$ be an edge of $\Gamma_X$ and identify $e$ with the unit
interval so that $\iota(e)\sim 0$ and $\tau(e)\sim 1$.  Let $p\colon
e\to X$ be a reduced edge path with $p(0)=\iota(e),$
$p^{-1}(\tau(e))=\emptyset$. The \term{fold of $X$ at $e$ along $p$}
is the space obtained by crushing any edges of $X'=X/(t\sim p(t))$ which
meet a valence one vertex of $\Gamma_{X'}$.


\begin{lemma}[\cite{stallings0}]
  \label{foldtoimmersion}
  If $\varphi\colon(X,\YY(X))\to(X',\YY(X'))$ is
  $\pi_1$--injective and combinatorial, elements of $\YY(X')$
  aspherical, then there is a relative graph
  $(\overline{X},\YY(\overline{X})),$
  $\YY(\overline{X})\equiv\YY(X),$ a homotopy
  equivalence $F\colon X\to\overline{X},$ and an immersion
  $\overline{\varphi}\colon\overline{X}\to X'$ such that
  $\overline{\varphi}\circ F$ is homotopic to $\varphi$ and $F$ is a
  composition of folds.
\end{lemma}


\par The proof of Lemma~\ref{foldtoimmersion} is an easy riff on
folding. Note that if $\varphi\colon X\to X'$ is $\pi_1$--injective
then it is homotopic to a combinatorial $\varphi'\colon
(X,\YY(X))\to(X',\YY(X'))$. This is accomplished by first homotoping
$\varphi$ so that $\YY_{\bullet}(X)$ has image in
$\YY_{\bullet}(X')$. This can be done since each $Y\in\YY(X)$ has
freely indecomposable fundamental group, $\varphi$ is
$\pi_1$--injective, and each $Y\in\YY(X')$ is aspherical. Then
homotope $\varphi$ so that every vertex of $\Gamma_X$ has image in
$\Gamma^{0}_{X'}$. Now subdivide the edges of $\Gamma_X$ and homotope
$\varphi$ to a combinatorial map.

\begin{proof}
  Suppose $\varphi$ isn't an immersion and satisfies condition 1 of
  Definition~\ref{immersion}. If $v$ maps to $w$ and
  $\pi_1(\varphi_w)$ isn't injective, then perform an ordinary
  Stallings fold.

  Suppose $\varphi$ doesn't satisfy condition 1, but does satisfy
  condition 2, and let $Y,$ $\st(Y),$ $N,$ and $V$ be as in
  Definition~\ref{immersion}. Let $b_1$ and $b_2$ be two vertices of
  $\partial N$ with a path $p\colon\left[0,1\right]\to N$ such that
  $p(0)=b_1$ and $p(1)=b_2,$ and such that $\left[\varphi\circ
    p\right]$ is trivial in $\pi_1(\st(Y),\partial\st(Y))$.

  Since $\varphi$ is $\pi_1$--injective, $b_1\neq b_2$. Adjacent to
  $b_1$ and $b_2$ are unique distinct edges $e_1,e_2\subset
  N,\Gamma_{X}$ so that $\iota(e_i)=b_i$. The path $p$ can be
  homotoped to the composition of two paths, the first traversing
  $e_1,$ and the second a path $p'\colon\left[0,1\right]\to
  N\setminus(e_1\setminus\tau(e_1))$ satisfying
  $p'^{-1}(\tau(e_1))=\set{0}$ and $p'(1)=\iota(e_2)=b_2$. Homotope
  the restriction of $\varphi$ to $e$ so it agrees with
  $(p')^{-1}$. Now let $X''$ be the fold of $X$ at $e_1$ along the
  path $p',$ with $F$ the quotient map. There is an obvious map
  $\varphi'\colon X''\to X'$ and the composition $\varphi'\circ F$ is
  homotopic to $\varphi$. Rinse and repeat. The process must terminate
  since the number of edges of $\Gamma_{X''}$ is strictly less than
  the number of edges of $\Gamma_{X}$.
\end{proof}

\par One important property of immersions of relative graphs is that
if $p$ is a reduced edge path then the composition of $p$ and an
immersion is also a reduced edge path. More generally, compositions of
immersions are immersions.

\subsection{Graphs of Spaces}

\label{graphsofspaces}

\subsubsection{Graphs of Free Groups}


\par In the next sub-subsection Theorem~\ref{maintheorem} will be
reduced to an analysis of spaces arising by adjoining roots to free
groups. We begin with a slightly more general construction than the
one we need since the analysis will give easy proofs of
Theorem~\ref{improvedshenitzer} and
Corollary~\ref{maximumrankconjugacy}. We now state the main theorem of
this subsection. In conjunction with Theorem~\ref{forthcomingsetup},
Theorem~\ref{addroottofreegroup} implies Theorem~\ref{maintheorem}

\begin{theorem}
  \label{addroottofreegroup}
  Let $X$ be a $2$--covered graph of spaces arising from adjunctions
  of roots to non-conjugate, indivisible elements $\gamma_i$ of a free
  group. Furthermore, suppose that $\gamma_i$ and $\gamma_j^{-1}$
  are nonconjugate for all $i\neq j$. If
  $\chi(\Gamma(X))=\chi(\underlying(X))$ then the edge spaces of $X$
  are trees.
\end{theorem}

\par Definitions follow.

\par An edge of the graph underlying a graph of spaces is denoted by a
lower case letter, and the space associated to that edge is denoted by
the same letter upper-cased. Edges of graphs of spaces are oriented,
and the edge map associated to the preferred orientation is typically
denoted ``$\tau$''.

\begin{definition}[Graph, Graph of Spaces]
  A graph is a set $W$ with an involution $\bar{\ }$ (``bar'') and
  retractions $\tau,\iota\colon W\to\mathrm{Fix}(\bar{\ })$
  compatible with $\bar{\ }$:
  \[\iota(\bar{w})=\tau(w),\quad \tau(\bar{w})=\iota(w)\]
The elements of $W\setminus\mathrm{Fix}(\bar{\ })$ are the
oriented edges of the graph. The fixed set of $\bar{\ }$ is the set of
vertices, and the maps $\tau$ and $\iota$ are the terminal and initial
vertices of oriented edges, respectively. We say that an edge $e$ is
\term{incident} to $v$ if $\tau(e)=v$.
\end{definition}

\par Note that a graph in this sense is a special kind of category. A
(ordinary) graph of spaces is a functor from a graph
$( W,\iota,\tau,\bar{\ })$ to $\mathrm{Top}$.

\par Let $\mathscr{G}$ be the category of simplicial graphs whose maps
are combinatorial immersions. For us, a graph of spaces is a functor
from a graph $( W,\tau,\iota,\bar{\ })$ to
$\mathscr{G}$. We'll be mostly interested in graphs of spaces which
satisfy a rather restrictive criterion on collections of edges
incident to vertices.

\par Members of a graph will be referred to with lower case variables,
and their images in $\mathscr{G}$ will have capital variable names.
If an edge $e$ is incident to $v,$ then we say that $E$ is
incident to $V,$ similarly for variables with subscripts.

\begin{definition}[$2$--Covered]
  A finite graph $V$ is 2-Covered by $\set{E_i}$ if, for every
  $i,$ there is an immersion $\tau_i\colon E_i\immerses V$ and each
  edge $f$ of $V$ is the image under $\coprod\tau_i$ of exactly
  two edges from $\coprod E_i$.  
\end{definition}

\par For the remainder of this section a graph of spaces $X$ will
satisfy the condition that if $e_1,\dotsc,e_n$ are incident to $v,$
then $V$ is 2-covered by $\set{\tau_i\colon E_i\immerses V}$.

\par We fix some notation for graphs of spaces.

\begin{itemize}
\item An underlying graph \underlying\ that the graph of spaces is
  built on, i.e., if $X$ is a graph of spaces, then $X$ really
  corresponds to a functor $\underlying(X)\to\mathscr{G}$.
\item vertex spaces are connected graphs $V_i$.
\item Edge spaces in the topological realization are products of
  intervals with connected ``edge-graphs'' $E_j$. Edge graphs may be
  points. Each edge space $E_j\times\mathrm{I}$ in the topological
  realization of $X$ has an embedded copy of $E_j,$
  $E_j\times\set{\frac{1}{2}}$.
\end{itemize}

\par For a graph of spaces, there is a natural (not necessarily
connected) subcomplex $\Gamma(X)$ consisting of horizontal edges:
$\Gamma(X)^{(0)}=X^{(0)},$
$\Gamma(X)^{(1)}=\bigcup_j E_j^{(0)}\times\mathrm{I},$ with
identifications induced by the immersions $\tau$. The
\term{horizontal} subgraph $\Gamma(X)$ is the realization of the
graph of spaces induced by restricting to the zero skeleta of the
vertex and edge spaces. (Note that we cheated a little. When we
defined a graph of spaces, we insisted on having connected vertex and
edge spaces. Zero skeleta are rarely connected, but the definition
makes sense just the same.) Let $\compcomp(X)$ be the subset of
$\Gamma(X)$ consisting of the connected components not
homeomorphic to $S^1$. Let $\sonecomp(X)$ be the subset
consisting of components homeomorphic to $S^1$.


\par For a graph of free groups over cyclic subgroups with a
homomorphism to a free group which embeds the vertex groups, there is
a natural complex $X$ which has a graph of spaces structure transverse
to the graph of spaces structure induced by its decomposition as a
graph of groups. This transverse graph of spaces structure is a
$2$--covered graph of spaces.

\par Let $ G=\Delta(F_i,Z_j)$ be a graph of free groups
$F_i$ (not necessarily nonabelian!)  over nontrivial cyclic
subgroups $Z_j$. If $\phi\colon G\onto\free_n$ is a homomorphism which
embeds each $F_i,$ then we can build a nice graph of spaces
representing $\phi$. For each $i,$ choose an immersion
$\varphi_i\colon\Gamma_i\immerses R_n,$ where $R_n$ is the rose with
$n$ petals, and fundamental group $\free_n$ and
$\pi_1(\Gamma_i)=F_i$. Each cyclic edge group $Z_j$ must embed in
$\free_n,$ so for each $j,$ choose an immersion
$\varphi_j\colon S^1_j\immerses R_n$ representing the image of
$Z_j$. If $Z_j\into F_i$ then $\varphi_j$ lifts to an immersion
$\varphi_{i,j}\colon S^1_j\immerses\Gamma_i$ (There may be more than
one possibility for $\varphi_{i,j},$ corresponding to a monogon in
$\Delta$. Choose two, one for each orientation of the edge.) Use the
data $\varphi_{i,j}$ to attach annuli, one for each edge of $\Delta,$
to the graphs $\Gamma_i,$ to build a graph of spaces $X$. Our
original homomorphism $\phi$ induces a map
$\varphi\colon X\to R_n$. Restricted to an annulus
$ S^1_j\times\mathrm{I},$ the map is projection to the first factor,
followed by the immersion $\varphi_i$. Let $b$ be the basepoint of
$R_n$. Now regard $X$ as a 2-covered graph of spaces by setting
$\set{V_p}$ to be the connected components of $\varphi^{-1}(b),$ and
edge graphs connected components of preimages of midpoints of edges of
$R_n$. The homomorphism $\phi\vert_{F_i}$ factors through the
inclusion $\Gamma_i\to X$. See the bottom two rows of
Figure~\ref{buildgraphofspaces}

\par The graph of groups $G$ has corank at most
\[\sum_{i=1}^l \rk(F_i)-l+1\] If $n=-\sum\chi(\Gamma_i)+1,$ then
$\chi(\underlying(X))=1-n=\sum\chi(\Gamma_i)=\chi(\Gamma(X))$.

\begin{lemma}
  \label{sillyrankinequality}
  If $ G=\Delta(F_i,Z_j)$ is a graph of free groups over nontrivial
  cyclic subgroups and $ G\onto\free_m$ then \[m\leq
  1-\sum\chi(F_i)\]
\end{lemma}

\par A homomorphism $\phi\colon\Delta(F_i,Z_j)\onto\free$ such that
the inequality of Lemma~\ref{sillyrankinequality} is an equality has
\term{maximal corank}.

\subsubsection{Non-free groups}
\label{nonfreegroups}


\par In this subsection we consider the situation of
Theorem~\ref{maintheorem}. The objective of this section is to prove
the following theorem.

\begin{theorem}
  \label{forthcomingsetup}
  Suppose that $\phi\colon G\into H$ and $H$ is a quotient of
  $G'=\adjoinroot{G}{\gamma_i}{k_i},$ $\boldsymbol{\gamma}_i,$ classes
  of indivisible elements of $G$ such that
  $\boldsymbol{\gamma}_i\neq\boldsymbol{\gamma}_j^{\pm 1}$ for all
  $i\neq j$ and $\gamma_i\in\boldsymbol{\gamma}_i$. Then
  $\scott(G)\geq\scott(H)$. If equality holds and $H$ has no $\zee_2$
  free factors then there is a partition of
  $\set{\boldsymbol{\gamma}_i}$ into subsets
  $\boldsymbol{\gamma}_{j,i},$ $j=0,\dotsc,p,$ $i=1,\dotsc,i_j,$ such
  that $\gamma_{j,i}\in G_j,\boldsymbol{\gamma}_{j,i},$ $j\geq 1,$ and
  a $2$--covered graph of spaces $X$ such that:
  \begin{itemize}
  \item $\chi(\underlying(X))=\chi(\Gamma(X))$.
  \item $\compcomp(X)$ is connected and has a collection of connected
    subgraphs $R_Y,$ $Y\in\YY(X_G),$ which have pairwise
    disjoint image under the map $\Gamma(X)\to\underlying(X)$.
  \item $H_j$ is a quotient of $\adjoinroot{G_j}{\gamma_{j,i}}{}$.
  \item For each $Y$ there is an attaching map $\psi_Y\colon R_Y\to
    Y$. The attaching maps of those components of $X$ with
    $\chi(\underlying)=0$ are the boundaries of mapping cylinders
    corresponding to $\gamma_{j,i},$ $j\geq 0$. The boundary of the
    mapping cylinder associated to $\gamma_{j,i}$ is attached to $Y_j$
    along $\gamma_{j,i}$.
  \item The space
    $X_G=(\compcomp(X)\sqcup\bigsqcup R_Y)/(x\sim\psi_Y(x))$
    has fundamental group $G$.
  \item The space
    $\widetilde{X}=(X\sqcup\bigsqcup R_Y)/(x\sim\psi_Y(x))$
    has fundamental group $G'$.
  \item The attaching maps $\gamma_{0,i}\to X_G$ factor through
  $\compcomp(X),$ and in fact $X$ is the union of $\compcomp(X)$ and
  the mapping cylinders for $\gamma_{0,i}$.
  \end{itemize}
\end{theorem}

\par The first and second bullets are key. We now prove
Theorem~\ref{maintheorem}, assuming Theorems~\ref{forthcomingsetup}
and~\ref{addroottofreegroup}.

\begin{proof}[Proof of Theorem~\ref{maintheorem}]
  The proof is by induction on $\scott(G)$. Assume the conclusions of
  Theorem~\ref{forthcomingsetup}. By Theorem~\ref{addroottofreegroup}
  every edge space of $X$ is a tree. Since the graphs $R_Y$ have
  disjoint image in the projection from $X$ to $\underlying(X)$ there
  is an edge $e$ of $\underlying(X)$ not in the image of any
  $R_Y$. The edge space associated to $E$ is a tree, hence there is an
  edge of $\compcomp(X)$ not attached to any $Y\in\YY(X_G)$
  and which is crossed exactly once by the representative of one of
  $\boldsymbol{\gamma}_{0,i},$ say $\gamma_{0,1}$. Then $G\cong
  G_1*\langle\gamma_{0,1}\rangle$, with all $\gamma_i\neq\gamma_{0,1}$
  conjugate into $G_1$. Let $G'_1=\adjoinroot{G_1}{\gamma_i}{k_i},$
  $\gamma_i\neq\gamma_{0,1}$ and let $H_1$ be the image of
  $G'_1$. Then $H\cong H_1*\langle
  \sqrt[k_{0,1}]{\gamma_{0,1}}\rangle$ and
  $\scott(G_1)=\scott(H_1)=(q_G-2,p_G)$. Repeating this procedure for
  all $\gamma_{0,i},$ reduce to the case where each $\gamma$ is
  conjugate into some $G_j$.

  Then $G_{i_0}=G_1*\dotsb*G_p*F$ and all leftover $\gamma_i$ are
  elements of some $\boldsymbol{\gamma}_{i,j},$ $i\geq 0,$ and
  \[G'_{i_0}=\adjoinroot{G_1}{\gamma_{1,i}}{}* \dotsb* 
  \adjoinroot{G_p}{\gamma_{p,i}}{}*F\] Passing to the image of
  $\adjoinroot{G_j}{\gamma_{j,i}}{}$ in $H,$ we see immediately that
  \[H_{i_0}\cong\img_H(\adjoinroot{G_1}{\gamma_{1,i}}{})*\dotsb*\img_H(\adjoinroot{G_p}{\gamma_{p,i}}{})*F\] Reassembling the free factors split off by Theorem~\ref{addroottofreegroup} proves the theorem.
\end{proof}

\begin{remark}
  Simply knowing that some $\gamma_{0,i}$ is a basis element in
  $\pi_1(\Gamma(X))$ isn't sufficient to imply that
  $\langle\gamma_{0,i}\rangle$ is a free factor of $G,$ thus
  Theorem~\ref{maintheorem} cannot be deduced from Baumslag's theorem.
\end{remark}

For the remainder of this section, fix
an inclusion $\phi\colon G\into H$ of finitely generated groups which
lifts to an epimorphism
$\widetilde{\phi}\colon G'\onto H$.  Before
we begin, replace $H$ by a group $H'$ as follows: Let $\phi\colon
L\onto L'$ be an epimorphism of groups. Then
\[\scott(\phi)=\max\set{\scott(L'')}:{\phi\mbox{ factors through an
epimorphism }L''\onto L'}\] Let $H'$ be a group such that $G'\onto H$
factors through $H'$ and such that $H'$ achieves $\scott(G'\onto
H)$. Clearly the freely indecomposable free factors of $H'$ have freely
indecomposable image in $H$.

\begin{lemma}
  \label{scottinequalityepi}
  Let $H'\onto H$ be as above. Then \[(q_{H'},p_{H'})\geq(q_H,p_H)\]

  If a freely indecomposable free factor of $H'$ has trivial image in
  $H$ then the inequality is strict.
\end{lemma}

\begin{proof}
  Let $\free_H$ be the free part of $H$. Then every freely
  indecomposable free factor of $H'$ is in the kernel of
  $H'\onto\free_H$ and we see that $\free_{H'}\onto\free_{H}$. Suppose
  that some freely indecomposable free factor of $H$ does not contain
  the image of a freely indecomposable free factor of $H'$. Let $H_0$
  be this freely indecomposable free factor. Then, by the reasoning
  above, $\free_{H'}$ maps onto $\free_{H}*H_0$. Since $\free_H$ has
  the same rank as $\free_{H'},$ $H_0$ must be trivial. Thus, $p_H\leq
  p_{H'}$.
\end{proof}

Rather than work with the inclusion $G\into H,$ we work with $G\into
H',$ suppressing the $'$ for convenience. 

Our first task is to find a suitable way to represent $\phi$ and
$\widetilde{\phi}$ as maps of cell complexes. We start by representing
$\phi$ as an immersion $\varphi\colon X_G\to X_H$ given by
Lemma~\ref{foldtoimmersion}. Once this is done, we build a nice space
$\widetilde{X}$ with fundamental group
$\adjoinroot{G}{\gamma_i}{k_i}$. This space is equipped with a well
behaved map to $X_H$ and we use this map to endow $\widetilde{X}$ with
a new graph of spaces structure transverse to the old one.

\par An admissible relative graph of spaces is \term{minimal} if it
has no valence one vertices, and for every valence two vertex $v$ of
$\Gamma_X\setminus\YY_{\bullet},$ the relative graph obtained
by unsubdividing $\Gamma_X$ at $v$ is inadmissible, and for every
valence one vertex $v$ of $\Gamma_X,$ $v$ is contained in some
$Y\in\YY(X)$.

\par
\begin{figure}[h]
  \psfrag{Y}{$Y$}
  \psfrag{Nepy}{$\st(Y)\in\mathcal{N}(X)$}
  \psfrag{einE}{$e\in\E(X)$}
  \psfrag{Nep}{$\st(v)$}
  \psfrag{v}{$v$}
  \psfrag{l}{$\{$}
  \psfrag{r}{$\}=X^I$}
  \psfrag{in}{}
  \centerline{%
    \includegraphics{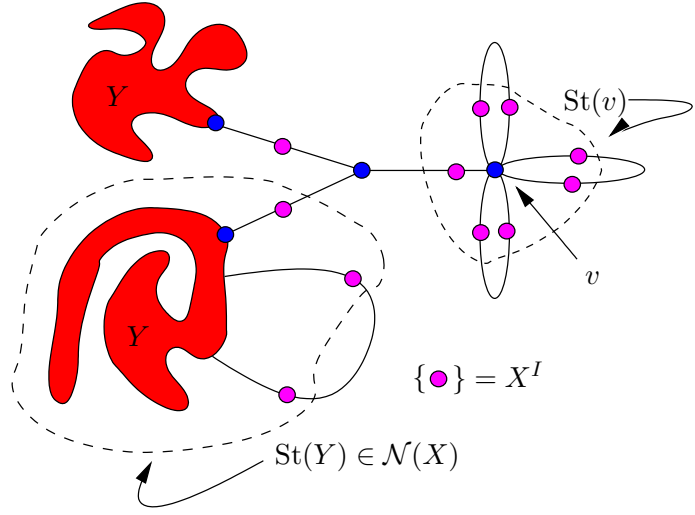}
  }
  \caption{A minimal admissible relative graph.}
  \label{normalform}
\end{figure}

 Suppose $X$ is admissible and minimal. Let
$\mathcal{N}(X)=\set{N_i}$ be the collection of all closures of
connected components of $X\setminus X^I$. If $N$ contains some element
$v\in X^E$ then $N$ is simply $\st(v)$. If $N\in\mathcal{N}(X)$ then
define $\partial N=N\cap X^I$.  A minimal admissible relative graph is
illustrated in Figure~\ref{normalform}.

\par Choose a relative graph of spaces $X_H$ with fundamental group
$H,$ and whose components $Y_1,\dotsc,Y_p\in\YY(H)$ are
aspherical and have fundamental group $\pi_1(Y_i)\cong H_i$. Choose a
relative graph $X_G$ with fundamental group $G$ and an immersion
$\varphi\colon X_G\immerses X_H$ provided by Lemma~\ref{foldtoimmersion}
representing $\phi$. We now build a space $\xtilde$ with fundamental
group $\adjoinroot{G}{\gamma_i}{k_i}$. Choose, for each $\gamma_i,$ an
immersion $\gamma_i\colon S^1\immerses X_G$ representing the conjugacy
class $\left[\gamma_i\right]$. If $\gamma$ is conjugate into some
$G_i$ then the immersion $\gamma\colon S^1\to X_G$ has image in the
connected component of $\YY_{\bullet}(X_G)$ representing
$G_i$. Let $\set{\gamma_i^j}$ be the subcollection of $\set{\gamma_i}$
consisting of elements conjugate into $G_i$. Let
$\set{\alpha_1,\dotsc,\alpha_m}$ be the subcollection of $\set{\gamma_i}$
which are not conjugate into any $G_i$.

\par Let $M_i$ be the mapping cylinder of the $k_i$ fold cover
$ S^1_i\to S^1,$ and let $r_i$ be the core curve (the range $ S^1$)
of $M_i$. Now glue the $M_i$ along $ S^1_i$ to $X_G$ using the
immersions $\gamma_i$ as attaching maps to form a space $\xtilde$. For
each $r_i,$ choose an immersion $r_i\to X_H$ representing the conjugacy
class of $\sqrt[k_i]{\gamma_i}\in H$. The map
$\varphi\colon X_G\immerses X_H$ lifts to a continuous map
$\widetilde{\varphi}\colon\xtilde\to X_H$ which agrees with the
immersions $r_i$ and $\varphi$.

\par Each mapping cylinder $M_i$ is a quotient space of an annulus
$A_i$. Let $\mathcal{A}$ be this collection of annuli and let
$i_a\colon\mathcal{A}_{\bullet}\to\xtilde$ be the disjoint union of
the maps $A_i\to M_i\to\xtilde$.


\par The space $\xtilde$ is a quotient of the disjoint union of of
$\mathcal{A}_{\bullet},$ $X_G,$ and the core curves $r_j$ of the
mapping cylinders $M_j$ (which are built out of annuli from
$\mathcal{A}$ and said core curves.). In analogy with a 2-covered
graph of spaces, set $\Gamma(\xtilde)=X_G\sqcup\bigsqcup
r_i$. Homotope $\widetilde{\varphi},$ relative to $\Gamma(\xtilde),$
so that $\mathcal{A}_{\bullet}\to X_H$ is transverse to $X_H^I$. Then
$i_a^{-1}(\mathcal{B}_{\bullet})=V$ is an embedded 1-submanifold of
$\mathcal{A}_{\bullet}$. Using an innermost disc argument and
asphericity of each component of $X_H^E,$ homotope
$\widetilde{\varphi}$ so that the submanifold contains no simple
closed curves bounding disks. If an (innermost) arc $\alpha\subset V$
has endpoints in only one boundary component of
$\mathcal{A}_{\bullet}$ then one of $X_G\to X_H$ or some $r\to X_H$
fails to be an immersion, contrary to hypothesis. If $V$ contains an $
S^1$ which does not bound a disk then some $\gamma_k$ vanishes in $H,$
contrary to hypothesis. Thus every connected component of $V$ is an
arc connecting distinct boundary components of some
$A\in\mathcal{A}$. If $\gamma_i$ is conjugate into some $ G_j$ then
the immersions $\varphi\circ\gamma_i$ and $r_i$ have images in
$\YY_{\bullet}(X_H)$ and any component of $V$ contained in $M_i$ is a
circle which does not bound a disk. A circle contained in an annulus
is homotopic to each boundary component, hence if such a circle
existed then $\gamma_i$ would have to be trivial in $G,$ a
contradiction. Thus, for such $M,$ $V\cap M=\emptyset$.


\par The first step in our analysis of $\widetilde{X}$ is to resolve
$\widetilde{\varphi},$ giving $\widetilde{X}$ a graph of spaces
structure transverse to the graph of spaces decomposition
$X_G\sqcup\bigsqcup M_i/\sim$. Let $\ZZ$ be the collection of
connected components of $\widetilde{\varphi}^{-1}(N),$
$N\in\mathcal{N}(X_H),$ and let $\mathcal{B}$ be the collection of
connected components of preimages of $X_H^I$. By transversality
$i_a^{-1}(\mathcal{B}_{\bullet})$ is a one manifold with boundary
contained in $\partial\mathcal{A}_{\bullet},$ and thus
$\mathcal{B}_{\bullet}$ is a graph contained in $\xtilde$. For each
$B\in\mathcal{B}$ there are two embeddings
$B\into\ZZ_{\bullet}\sqcup\mathcal{W}_{\bullet}$. Note that if
$\gamma_i$ is conjugate into some $ G_j$ then $M_i$ is completely
contained in some $Z$. The boundary of $Z\in\ZZ$ is the set
$\mathcal{B}_{\bullet}\cap Z,$ and coincides with the set of points of
$Z$ mapping to $X_H^I$. Inclusions $B\into Z$ are simply inclusions of
boundary components.

\begin{figure}[tb]
  \psfrag{scottinequality}{Lemma~\ref{scottinequality}}
  \psfrag{or}{or}
  \psfrag{v}{}
  \psfrag{w}{}
  \centerline{%
    \includegraphics[scale=0.5]{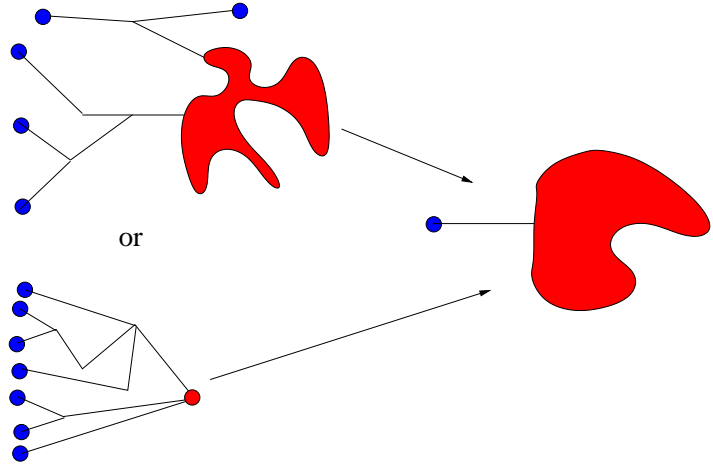}
  }
  \caption{Possibilities for $S$ mapping to
    $N\in\mathcal{N}(X_H)$.}
  \label{piecesmappingtorelativespaces}
\end{figure}

\par Some possibilities for $Z\cap X_G$ are illustrated in
Figure~\ref{piecesmappingtorelativespaces}.

\par Each mapping cylinder $M_j$ is either completely contained in
$\ZZ_{\bullet}$ or has nontrivial intersection with
$\mathcal{B}_{\bullet}$. If $M$ has nontrivial intersection with
$\mathcal{B}_{\bullet}$ then $r\cap Z$ (recall that $r$ is the core
curve of $M$) is a collection of closed intervals, even in number. The
preimage $i_a^{-1}(\mathcal{B}_{\bullet})$ slices an annulus
$A\in\mathcal{A}$ into rectangles
$R^z_1,R^e_1\dotsc,R^z_m,R^e_m,$ where the $R^z_*$ have images in
$\mathcal{N}_{\bullet}(X_H)$ and $R^e_*$ have images in
$\E_{\bullet}(X_H)$.

\par The connected components of the intersections of the relative
graph $\Gamma(\xtilde)$ and elements of $\ZZ$ or $\mathcal{W}$
are relative graphs with distinguished valence one vertices
$\set{Z\vert W}\cap\mathcal{B}_{\bullet}$. Let $\mathcal{S}(Z)$ be the
collection of connected components of $\Gamma(\xtilde)\cap Z$.  If
$S\in\mathcal{S}$ set $\partial S=S\cap\mathcal{B}_{\bullet}$.  For a
fixed $Z\in\ZZ$ let $\RR(Z)$ be the subcollection of
all $R^z_*$ (coming from all annuli $A\in\mathcal{A}$) which are
contained in $Z$.

The boundary of each rectangle $R$ is composed of two types of arcs,
$\partial^{\pm}R=R\cap\partial\mathcal{A}_{\bullet}$ and
$\partial^{L|R}R=R\cap i_a^{-1}(\mathcal{B}_{\bullet})$. The former
shall be known as horizontal boundary arcs and the latter as vertical
boundary arcs.

Let $\partial^+R$ be a horizontal boundary arc of some
$R\in\RR(Z)$. Then $\varphi^+=i_a\vert_{\partial^+ R}$ has
image in some connected component $S$ of $\mathcal{S}_{\bullet}(Z)$
and $\partial\partial^+ R$ maps to $\partial
S=S\cap\mathcal{B}_{\bullet}\subset\partial Z$. We define $\varphi^-$
similarly. If $S$ is a tree then this path is an embedding and
connects distinct boundary components of $S$. If $S$ contains a
relative space $Y$ of $X_G$ then it is a reduced edge path since the
maps representing $\alpha_i$ were immersions.


We reconstruct $Z$ by gluing the rectangles $\RR(Z)$ to
$\mathcal{S}_{\bullet}(Z)$ via the attaching maps
$\varphi^{\pm}\colon(\partial^{\pm}R,\partial\partial^{\pm}R)\to
\bigsqcup(S,\partial S)\in\mathcal{S}(Z)$.

\begin{figure}[h]
  \psfrag{psib}{$\psi^{-1}(b)=B\in\mathcal{B}$}
  \psfrag{psie}{$\psi^{-1}(e)$}
  \psfrag{Z}{$Z$}
  \psfrag{OZ}{$\zbar$}
  \psfrag{bbar}{$\overline{B}$}
  \centerline{%
    \includegraphics{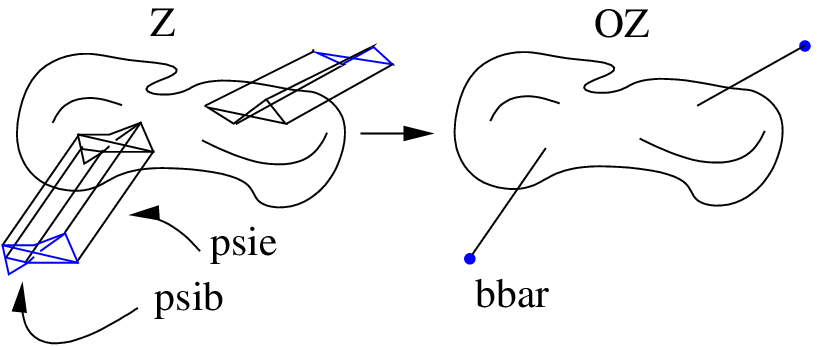}
  }
\caption{}
\label{xhprime}
\end{figure}

\par The boundary of $Z$ is the union of vertical boundary arcs of the
rectangles comprising it, along with all valence one vertices of
$\partial\mathcal{S}_{\bullet}(Z)$ not contained in some vertical
boundary arc of a rectangle. By construction, $Z$ is a connected
component of $\varphi^{-1}(N)$ for some
$N\in\mathcal{N}(X_H)$. Suppose $N$ is a star of some $v\in
X_H^E$. Let $\set{b_1,\dotsc,b_k}$ be the valence one vertices
comprising $\partial N,$ with incident edges $e_i\subset N$ so that
$\tau(e_i)=b_i$. If $\psi\colon Z\to N$ is the restriction of
$\widetilde{\varphi}$ then $\psi^{-1}(e_i)$ is a collar neighborhood
of $\psi^{-1}(b_i)\in\mathcal{B},$ a boundary component of $Z$. The
restriction $\psi$ factors through the map which projects each collar
onto the $\mathrm{I}$ factor. Call the resulting quotient space
$\zbar$. If $N$ contains elements of $X_H^E$ then $Z$ is homeomorphic
to a product $B\times\mathrm{I},$ for some $B$ in $\mathcal{B},$ let
$Z\to\overline{Z}$ be the projection to the $\mathrm{I}$ factor. For
each $B$ in $\mathcal{B}$ let $\overline{B}$ be the quotient space
consisting of a single point.

\par
\begin{figure}[h]
  \centerline{%
    \xymatrix{%
            \coprod B^{(0)}\ar[r]\ar[d] & \coprod\mathcal{S}_{\bullet}(Z)\ar[r]\ar[d] & X_G\sqcup\coprod r_i \ar[d] & \\
	    \coprod B\ar[r]\ar[d] & \coprod Z\ar[r]\ar[d] & \xtilde\ar[dr]\ar[d] & \\
	    \coprod \overline{B}\ar[r] & \coprod \zbar\ar[r] & X_K\ar[r] & X_H
    }
  }
  \caption{Construction of $X_K$. Horizontal diagrams are pushouts.}
  \label{constructxk}
\end{figure}
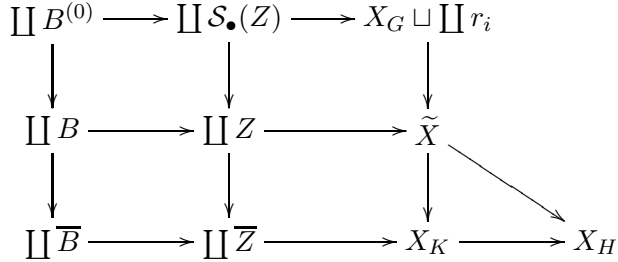

 The lift $\widetilde{\varphi}\colon\xtilde\to X_H=X_H$ factors
through the graph of spaces $X_K,$ $\pi_1(X_K)=K,$ obtained by
reassembling the collection $\set{\zbar}:{Z\in\ZZ}$ If
$B\subset Z_1,Z_2,$ then identify the images $\overline{B}\in\zbar_1$
and $\overline{B}\in\zbar_2$. Construction of $X_K$ is illustrated in
Figure~\ref{constructxk}.

\par Let $S_1,\dotsc,S_n$ be the connected components of
$\mathcal{S}(Z),$ and let $R_1,\dotsc,R_k$ be the rectangles from
$\RR(Z)$. For each rectangle let
$\varphi_j^{\pm}\colon\partial^{\pm}R_j\to\mathcal{S}_{\bullet}(Z)$ be
the attaching maps for the horizontal boundary arcs of $R_j$. Then $Z$
is the quotient space \[\mathcal{S}_{\bullet}(Z) \sqcup
\RR_{\bullet}(Z)/ (x\sim\varphi_j^{\pm}(x))\] The boundary of
$Z$ consists is the union $\bigcup\partial_j^{L\vert
R}R_j\cup\bigcup_i\partial S_i=\mathcal{B}_{\bullet}\cap Z\subset Z$.

\begin{lemma}
  The spaces $\zbar$ have freely indecomposable (or trivial!)
  fundamental groups.
\end{lemma}

\begin{proof}
  Suppose not. Then $\pi_1(\xtilde)\onto H$ factors through a group
  with strictly higher Scott complexity. Recall that $H$ is the new
  name for $H',$ chosen to achieve $\scott(G'\onto H)$. If some
  $\zbar$ had freely decomposable fundamental group then
  $\scott(K)>\scott(H')$.
\end{proof}


\par Recall that $q_H\leq q_G,$ and that if equality holds then no
$S\in\mathcal{S}(Z)$ has fundamental group with nontrivial free part,
thus we may assume that no $S$ has $\zee$ as a free factor of its
fundamental group.

\par Let $G=G(Z)$ be the graph with vertex set $\mathcal{S}(Z)$ and
edge set $\RR(Z)$. The endpoints of an edge $R$ are the
boundary components $\partial^{\pm}R,$ and an endpoint
$\partial^{\pm}R$ is attached to $S$ if the image of $\varphi^{\pm}$
is contained in $S$. Let $T$ be a maximal tree in $G$. Build a space
$Z_T$ by restricting to the tree $T$.

\par Since no $S\in\mathcal{S}(Z)$ has nontrivial free part and $G$
embeds in $H,$ the components $Z\in\ZZ$ fall into three
classes:

\begin{enumerate}
  \item $\pi_1(\zbar)$ is trivial. Such $Z$ contain no $Y\in\YY(X)$.
  \item $\pi_1(\zbar)$ is nontrivial and $Z$ contains no $Y\in\YY(X)$
  \item $\pi_1(\zbar)$ is nontrivial and $Z$ contains some $Y\in\YY(X)$
\end{enumerate}

\par Let $\ZZ_i$ be the subset of $\ZZ$ containing all
$Z$ of the $i$--th type. For each $Z,$ let
\[\Delta_q^-(Z)=\frac{1}{2}\left(\#\partial Z_T-\#\partial\zbar\right)\] 
and for $Z\in\ZZ_3$ let
\[\Delta_q^+(Z)=\frac{1}{2}\betti(\partial Z_T)\] and for
 $Z\in\ZZ_{1|2}$ set $\Delta_q^+(Z)=0$.

If $S$ is a relative graph with no loops (no contribution to $q_G$),
set $\kappa(S,\partial S)=\frac{1}{2}\#\partial S -1,$ and observe
that
\[\sum_{Z\in\ZZ}\sum_{S\in\mathcal{S}(Z)}\kappa(S,\partial
S)=q_G-1\] and
\[\sum_{Z\in\ZZ}\kappa(\zbar,\partial\zbar)=q_H-1\] The 
complexity $\kappa$ is intended to be a stand-in for curvature. Beware
the sign convention we've chosen. 

If $Z\in\ZZ_1$ then define $\Delta_p^+(Z)=\Delta_p^-(Z)=0$. If
$Z\in\ZZ_2$ then define $\Delta_p^+(Z)=1$ and
$\Delta_p^-(Z)=0$. If $Z\in\ZZ_3$ and $Y_1,\dotsc,Y_{k+1}$ are
the components of $\YY(X_G)$ contained in $Z,$ then
$\Delta_p^-(Z)=k$ and $\Delta_p^+(Z)=0$. We now give three lemmas
relating the quantities $\Delta_{p|q}^{\pm}$ to one another.

\begin{figure}[ht]
\psfrag{XFIGLABEL}{LaTeXSTUFF}
\psfrag{R}{$R$}
\psfrag{R-}{$\partial^-R$}
\psfrag{R+}{$\partial^+R$}
\psfrag{Rl}{$\partial^LR$}
\psfrag{Rr}{$\partial^RR$}
\psfrag{S}{$S$}
\centerline{%
  \includegraphics{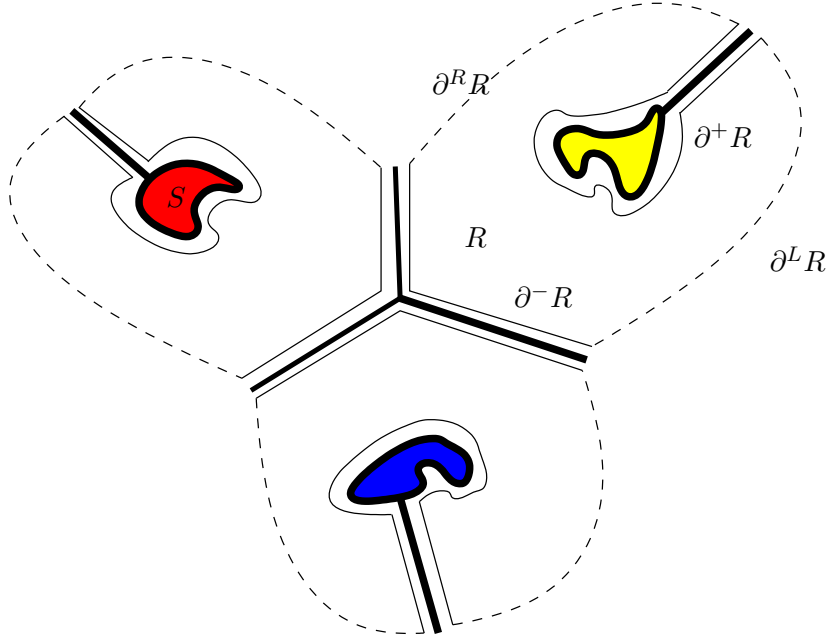}
}
\caption{Illustration for Lemma~\ref{deltapforz3}. Each independent
  loop in $\partial Z_T$ must contribute \emph{at least} $1$ to
  $\Delta_p^-$. The tree in this example is a tripod, there is one
  boundary component of $Z_T,$ and $\Delta_p^-=2$. A similar example
  with only two components $S$ would have $\Delta_p^-=1$.}
\label{deltapforz3figure}
\end{figure}

\begin{lemma}
  \label{deltapforz3}
  \[\kappa(\zbar,\partial\zbar)-\sum_{S\in\mathcal{S}(Z)}\kappa(S,\partial S)
  =\Delta_q^+(Z)-\Delta_q^-(Z)\]

  If $Z\in\ZZ_3$ then
  $\Delta_p^-(Z)\geq \betti(\partial Z_T) = 2\Delta_q^+(Z)$.
\end{lemma}

\begin{proof}
  To show the first equality we only need establish that
  \[\frac{1}{2}\#\partial Z_T -1 = 
  \sum_{S\in\mathcal{S}(Z)}(\frac{1}{2}\#\partial S -1) +\Delta_q^+(Z)\]

  To this end, let $T\subset G(Z)$ be a maximal tree, and let
  $S_1,\dotsc,S_k$ be an enumeration of $\mathcal{S}(Z)$ such that
  $S_{i+1}$ is connected to $S_0\cup R_0\cup S_1\cup R_1\dotsb S_i$ by an
  edge $R_i\subset T$. Assume that $R_i$ is oriented so that
  $\partial^+R_i$ is attached to $S_{i+1}$. Let $Z_i$ be the union of
  $S_1,\dotsc,S_i$ and $R_1,\dotsc,R_{i-1}$.

  The boundary of $\partial^{+|-}R_i$ consists of two points. Suppose
  that for at least one of $+$ or $-,$ the image
  $\partial\partial^{+|-}R_i$ is contained in two \emph{distinct}
  boundary components of at least one of $Z_i$ or $S_{i+1}$. If this
  is the case then
  \begin{align*}
  \kappa(Z_{i+1},\partial Z_{i+1}) &= \frac{1}{2}\#\partial Z_{i+1}-1 \\
                                    &= \frac{1}{2}\left(\#\partial Z_i + \#\partial S_{i+1}-2\right)-1\\ 
                                    &= \kappa(Z_i,\partial Z_i)+\kappa(S_{i+1},\partial S_{i+1})
  \end{align*}

  If both $\partial\partial^+R_i$ and $\partial\partial^-R_i$ have
  image in the same boundary component of $S_{i+1},$ $Z_i,$
  respectively, then $Z_i$ contains at least $\betti(\partial Z_i)+1$
  elements of $\YY$ and $S_{i+1}$ must contain a new element
  of $\YY$. We then have
  \begin{align*}
    \kappa(Z_{i+1},\partial Z_{i+1}) &= \frac{1}{2}\#\partial Z_{i+1}-1 \\
                                      &= \frac{1}{2}\left(\#\partial Z_i + \#\partial S_{i+1}-1\right)-1 \\
                                      &= \kappa(Z_i,\partial Z_i)+\kappa(S_{i+1},\partial S_{i+1}) + \frac{1}{2}
  \end{align*}
  Each such rectangle makes a contribution of $+1$ to $\betti(\partial
  Z_T),$ $+1/2$ to $\Delta_q^+,$ and only such rectangles make such
  contributions, thus
  \[\kappa(Z_T,\partial Z_T)=\sum_{S\in\mathcal{S}(Z)}\kappa(S,\partial S)+\Delta_q^+(Z)\]

  We now need to compare $\Delta_q^+(Z)$ to $\Delta_p^-(Z)$. If
  $\partial\partial^{\pm}R$ maps to a single boundary component of
  $S\in\mathcal{S}(Z)$ then, since $\partial^{\pm}R\to S$ is an
  immersion, $S$ must have nontrivial fundamental group, and since
  the free part of $S$ is trivial, it must contain some element $Y$ of
  $\YY(X_G)$. 

  Choose the exhaustion of $T$ so that $S_0$ has an incident edge $R,$
  $\partial^-R\to S_0$ such that $\partial\partial^-R$ maps to a
  single boundary vertex of $S_0$. By the reasoning above, $S_1$
  contains some element $Y_0\in\YY(X_G)$. Let
  $R_{i_1},\dotsc,R_{i_{\betti(\partial X_T)}}$ be the rectangles such
  that $\partial\partial^+R_{i_j}$ maps to a single boundary component
  of $S_{i_J+1}$. Then each $S_{i_j+1}$ contains some element
  $Y_{j}\in\YY(X_G)$. Since the $S_{i_j}$ are distinct, $Z$
  contains at least $\betti(\partial Z_T)$ elements of
  $\YY(X_G),$ i.e., $\Delta_p^-(Z)\geq\betti(\partial Z_T)$.
\end{proof}

\par Not all spaces which abstractly resemble $Z$'s occur as $Z$'s. We
now give a definition for a certain class of useful spaces resembling
them.

\begin{definition}
  \label{unionoftrees}
  A \term{union of trees} is a graph of spaces $Z$ whose vertex
  spaces are relative trees $S_1,\dotsc,S_n,$ edge spaces are
  intervals $I_1,\dotsc,I_m,$ and whose attaching maps $(I,\partial
  I)\to(S,\partial S)$ are reduced edge paths. The boundary of $Z$ is
  the union \[\partial Z=\bigcup\partial S_i\cup\bigcup(\partial
  I_j\times I)\]
\end{definition}

\begin{lemma}
  \label{deltaqforz2}
  \label{structurefortrees}
  If $\Delta_q^-(Z)=0$ then $Z$ has the following form: There are
  subcollections $\RR_i\subset\RR(Z)$ such that the
  restriction $Z_i$ of $Z$ to the rectangles $\RR_i$ is
  homeomorphic to a product $G_i\times I,$ $G_i$ a graph, and $Z$ is
  recovered by gluing components of $\sqcup(G_i^{(0)}\times I)$ to
  $\mathcal{S}_{\bullet}(Z)$. The graph with vertex set
  $\set{Z_i}\cup\mathcal{S}(Z)$ and an edge between $Z_i$ and $S$ if
  $v\times I\subset Z_i$ is identified with an edge path in $S$ is a
  tree.

  If $Z\in\ZZ_2$ then $\Delta_q^-(Z)\geq\frac{1}{2}$. If
  equality holds then $\pi_1(\zbar)\cong\zee_2$.
\end{lemma}

\begin{proof}
  First suppose that $\Delta_q^-(Z)=0$. As before, let $G(Z)$ be the
  graph with vertex set $\mathcal{S}(Z),$ edge set $\RR(Z),$
  and maximal tree $T\subset G$. Let $Z_T$ be the graph of spaces
  obtained by restricting to $T$ and let $\RR'(Z)$ be the
  subset of $\RR(Z)$ consisting of rectangles not contained in
  $T$.

  The boundary of $Z_T$ is a forest in the boundary of $Z$. Let $\sim$
  be the equivalence relation on $\RR$ generated by $R_1\sim R_2$
  if \[ \left[\partial^{+}R_1\to\mathcal{S}(Z)\right] =
  \left[\partial^{+}R_2\to\mathcal{S}(Z)\right]\] Let $\sim'$ be the
  same equivalence relation restricted to $Z_T$. For each
  $\sim'$--equivalence class $\left[R\right]$ let $Z_{\left[R\right]}$
  be the subspace of $Z_T$ obtained by restricting to rectangles in
  $\left[R\right]$. Let $B$ be a boundary component of
  $Z_{\left[R\right]}$. Then $Z_{\left[R\right]}$ is homeomorphic to
  the product $B\times I$.

  We claim that the map from $\sim'$ equivalence classes to $\sim$
  equivalence classes is an injection. Consider a rectangle
  $R\in\RR',$ let $S^+$ ($S^-$) be the member of $\mathcal{S}$
  containing the image of $\partial^+R$ ($\partial^-R$), and let
  $R_1,\dotsc,R_n$ be the path in $T$ from $S^-$ to $S^+$. The
  configuration of $R_i$ and $S^{\pm}$ has the form illustrated in the
  following figure. Boundary components are bold.

  \begin{figure}[h]
    \psfrag{R1}{$R_1$}
    \psfrag{dots}{$\dotsb$}
    \psfrag{Rn}{$R_n$}
    \psfrag{S+}{$S^+$}
    \psfrag{S-}{$S^-$}
    \psfrag{}{}
    \psfrag{}{}
    \psfrag{}{}
    \centerline{%
      \includegraphics{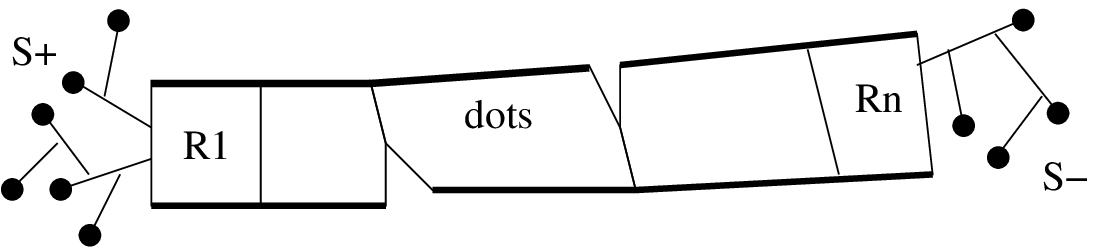}
    }
  \end{figure}
  
  The boundary components of $S^{\pm}\cup R_i$ are contained in
  \emph{disjoint} boundary components of $Z'$. Attaching $R$ to this
  configuration, we see that to satisfy $\Delta_q^-(Z)=0,$
  $\partial^+R$ must be attached to $\partial^+R_n,$ $\partial^-R$
  must be attached to $\partial^-R_1,$ for every $i,$ $R_i\sim
  R_{i+1},$ and in fact the induced orientations of $\partial^{\pm}R$
  must be coherent. From this it is easy to see that $\sim$
  equivalence classes must be homeomorphic to the products $B\times I$
  above. To construct $\zbar,$ collar neighborhoods of boundary
  components of $Z$ are crushed to intervals. From the
  characterization of $\sim$ equivalence classes above, we see that
  $\zbar$ is homotopy equivalent to a wedge of spheres, hence has
  trivial fundamental group.

  The decomposition of $Z$ as a union of products follows immediately.

  Now suppose that $\Delta_q^-(Z)>0$. Clearly
  $\Delta_q^-(Z)\geq\frac{1}{2}$. Suppose $\Delta_q^-(Z)$ is
  $\frac{1}{2}$. Let $Z_T$ be as before, and attach a rectangle $R$ to
  $Z_T$ such that $\#\partial Z_T\cup R=\#\partial Z_T-1$. Let
  $R_1,\dotsc,R_n,$ $S^{\pm}$ be a path in $T$ as before. Then $R$ and
  $R_i$ must be in the configuration illustrated in
  Figure~\ref{rpath3}. The boundary is in bold. Adding rectangles,
  maintaining $\Delta_q^-(Z)=\frac{1}{2},$ does not change the
  fundamental group of $\zbar,$ which is clearly $\zee_2$. In
  particular, \[\sum_{Z\in\ZZ_2}\Delta_q^-\geq\Vert\ZZ_2\Vert\]

  \begin{figure}[h]
    \psfrag{R1}{$R_1$}
    \psfrag{R2}{$R_2$}
    \psfrag{Rn-1}{$R_{n-1}$}
    \psfrag{Rn}{$R_n$}
    \psfrag{R}{$R$}
    \psfrag{insplus}{$\subset S^-$}
    \psfrag{insminus}{$\subset S^+$}
    \centerline{%
      \includegraphics[scale=0.5]{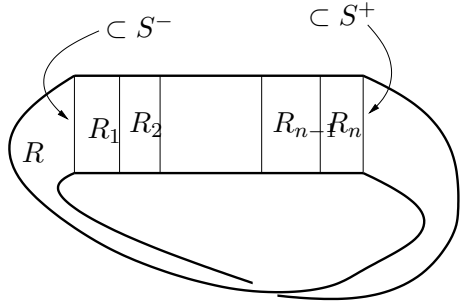}
    }
    \caption{Illustration for Lemma~\ref{deltaqforz2}.}
    \label{rpath3}
  \end{figure}

\end{proof}

Since $q_H\leq q_G$ we see immediately that if $q_H=q_G$ then
\[\sum_{Z\in\ZZ_3}\Delta_q^+(Z) =
\sum_{Z\in\ZZ}\Delta_q^-(Z)\] Then there are \emph{at most}
$2\sum_{Z\in\ZZ_3}\Delta_q^+(Z)$ elements $Z\in\ZZ_2$
since each such component has $\Delta_q^-(Z)\geq\frac{1}{2}$. Thus
$\Vert\ZZ_2\Vert=\sum\Delta_p^-(Z)\leq\sum_{Z\in\ZZ_3}\Delta_p^-(Z)$
and we conclude that $p_{H'}\leq p_G$. If equality holds then for each
$Z\in\ZZ_2,$ $\Delta_q^-(Z)=\frac{1}{2},$ and for each
$Z\in\ZZ_3,$ $\Delta_p^-(Z)=1$. Moreover, by
Lemma~\ref{deltaqforz2}, each element $\overline{Z}$ has fundamental
group $\zee_2$ for $Z\in\ZZ_2,$ and we have the following
lemma.

\begin{lemma}
  \label{scottlemma}
  \label{scottinequality}
  If $q_{H'}=q_G$ then $p_{H'}\leq p_G$. For all $Z\in\ZZ_2,$
  $\Delta_q^-(Z)=\frac{1}{2}$. In particular, every such $\zbar$ has
  fundamental group $\zee_2$.
\end{lemma}





\par We now revert to $H=H$ and $H'=H'$. 

\begin{lemma}
  If $H$ has no $\zee_2$ free factors and $(q_G,p_G)=(q_H,p_H)$ then
  $\ZZ_2$ is empty, all members of $\ZZ_3$ contain
  exactly one element of $\YY(X_G),$ and $\Delta_{p\vert q}^{\pm}(Z)=0$
  for all $Z$.
\end{lemma}

\begin{proof}
  By Lemma~\ref{scottinequalityepi} $\scott(H)\leq\scott(H'),$ and by
  Lemma~\ref{scottlemma} $\scott(G)\geq\scott(H')$. Suppose
  $\ZZ_2$ is nonempty. Then $H'$ has a freely indecomposable
  free factor $\zee_2$. Since $H$ has no $\zee_2$ free factors, this
  implies $\scott(H)<\scott(G),$ contrary to hypothesis.

  Thus $\Delta_q^\pm(Z)=0$ for all $Z$. Let $Z\in\ZZ_3$. Since
  $\zbar$ is freely indecomposable and $\pi_1(\zbar)$ maps to a freely
  indecomposable free factor of $H,$ to have equality in the second
  coordinate, each such $Z\in\ZZ_3$ must contain exactly one
  $Y\in\YY(X_G)$.
\end{proof}


\begin{definition}
  An union of trees is \term{treelike} if $\Delta_q^-(Z)=0,$ as is an
  element of $\ZZ_3$ if $\Delta_q^-(Z)=\Delta_q^+(Z)=0$.
\end{definition}

Let $Z$ be a treelike union of trees, and express $Z$ as a union of
$\sqcup (B_i\times I)\cup\mathcal{S}(Z)$ modulo attaching maps. For
each $B\times I,$ let $\pi_B$ be the projection onto the $I$
coordinate. Give $B\times I$ the foliation whose leaves are
$\pi^{-1}(x),$ $x\in I,$ and give each $S\in\mathcal{S}(Z)$ the
foliation whose leaves are simply the points of $S$. Then define
$\mathcal{F}(Z)$ to be the foliation on $Z$ induced by the foliations
on $B\times I$ and $S$. Define $\Gamma_Z$ to be the leaf space of
$\mathcal{F}(Z),$ and denote the quotient map by $\pi_Z$. The
following lemma is obvious from the construction.

\begin{lemma}
  \label{extensionlemma}
  The following facts about $\Gamma_Z$ are true:
  \begin{itemize}
  \item $\Gamma_Z$ is a finite tree.
  \item Each $S\in\mathcal{S}(Z)$ embeds in $\Gamma_Z$ under the
    quotient map. Any two images intersect in at most an interval.
  \item Each valence one vertex of $\Gamma_Z$ is the image of exactly
    one boundary component of $Z$.
  \item $\kappa(\Gamma_Z,\partial\Gamma_Z)=\kappa(Z,\partial Z)$
  \item Point preimages under $\pi_Z$ are connected.
  \end{itemize}

  The following extension property also holds. Let $f\colon Z\to A$ be
  a continuous map to an aspherical space $A$. If every boundary
  component of $Z$ is mapped to a point and $g\colon(I,\partial
  I)\to\mathcal{S}_{\bullet}$ is a reduced edge path, then there
  exists a lift $\widetilde{f}\colon\Gamma_Z\to A$ such that
  $\widetilde{f}\circ\pi_Z$ is homotopic to $f$ via a homotopy which
  is constant on $\partial Z$ and $I$.
\end{lemma}

The space $\xtilde$ was constructed by adjoining mapping cylinders $M$
along their boundaries to elements represented by immersions
$\gamma\colon S^1\immerses X_G$. Any interesting homotopies of
$\gamma$ are supported on arcs $I=\left[a,b\right]\subset S^1$ such
that $\gamma(I)\subset\YY_{\bullet}(X_G)$. In fact, up to such
homotopies, $\gamma$ is essentially unique: if an immersion $\gamma'$
is chosen, rather than $\gamma,$ to represent the conjugacy class of
$\gamma\in G$ (again, conflating immersions and conjugacy classes),
the space $\xtilde'$ constructed differs from $\xtilde$ only in that
for some (possibly more than one) $Z\in\ZZ,$ the attaching map
of a rectangle $\partial^{\pm}R\to S$ is altered by a homotopy
supported on an arc contained in $\partial^{\pm}R$ and having image in
some $Y\in\YY$ (contained in $S$).

\begin{convention}
  \label{chooserepresentatives}
  For convenience, we choose, for every homotopy class $\left[p\right]$
  in $\pi_1(S,\partial S),$ $S\in\mathcal{S}(Z),$ $Z\in\ZZ,$ a
  unique representative reduced edge path $p\colon(I,\partial
  I)\to(S,\partial S),$ such that $\left[p\right]^{-1}$ is represented
  by $t\mapsto p(1-t)$ unless $\left[p\right]$ represents a two-torsion
  element, in which case we choose $p$ to represent
  $\left[p\right]^{-1}$. Construct $\xtilde$ so that every attaching map
  $\partial^{\pm}R\to\mathcal{S}_{\bullet}$ agrees with the chosen
  representative in its homotopy class.
  
  For each $Z,$ under the hypothesis that $(q_H,p_H)=(q_G,p_G),$ there
  is at most one element $S_0\in\mathcal{S}(Z)$ containing an element
  $Y\in\YY(X_G)$. Without loss of generality, we may assume
  that $S_0$ is the star of $Y$ in $X_G$ and $\Gamma_{X_G}\cap S_0$ is
  a single point. See Figure~\ref{yfigure}. For such $S_0$ let $e_i$
  be the oriented edges of $\Gamma_{X_G}$ such that
  $\tau(e_i)=\tau(e_j)=b\in Y$ for all $i,j$. Then $S_0$ takes the
  form
  \begin{figure}[h]
    \psfrag{e1}{$e_1$}
    \psfrag{e2}{$e_2$}
    \psfrag{en}{$e_n$}
    \psfrag{dots}{$\vdots$}
    \psfrag{Y}{$\ni Y$}
    \centerline{%
      \includegraphics{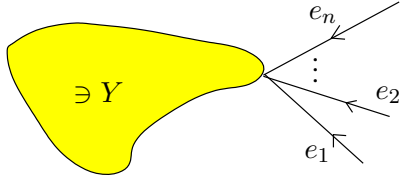}
    }
  \caption{Illustration for Convention~\ref{chooserepresentatives}.}
  \label{yfigure}
  \end{figure}
\end{convention}

\begin{lemma}
\label{zscontainingoney}
Suppose $Z$ is treelike and contains only one element $Y$ of
$\YY(X_G)$. Let $S_0$ be the element of $\mathcal{S}$
containing $Y$. Then there are treelike unions of trees
$Z_1,\dotsc,Z_n,$ $S_i\in\mathcal{S}(Z_i),$ $S_i\cong I,$ and reduced
edge paths $h_i\colon (S_i,\partial S_i)\to (S_0,\partial S_0)$ (which
are in the fixed list of representatives of homotopy classes
$\partial^{\pm}R\to\mathcal{S}_{\bullet}$) such that
\[Z=((S_0\sqcup\bigsqcup Z_i)/(x\sim h_i(x)))\]
\end{lemma}

The proof of Lemma~\ref{zscontainingoney} will resemble the proof of
Lemma~\ref{structurefortrees}.

\begin{proof}
  As in Lemma~\ref{structurefortrees}, construct $G(Z)$. Let
  $G_1,\dotsc,G_n$ be the closures of the connected components of
  $G\setminus\set{S_0},$ and let $T_i$ be a maximal tree in $G_i$ which
  meets $S_0$ only once. Construct $Z_i'$ by restricting $Z$ to $T_i,$
  an consider what happens when a rectangle is attached to $Z_i'$ (it
  must be attached along both horizontal boundary arcs to $Z_i'$).

  If $R$ is not attached to $S_0,$ we carry out the same analysis done
  in Lemma~\ref{structurefortrees}, which we now revisit. Consider
  $R\in G_i\setminus T_i$ which is attached along $\partial^-R$ to
  $S_0,$ and construct the path $R_1,\dotsc,R_n$ from $S^-$ to $S^+,$
  and assume that $S^-=S_0$. The union of the $R_i$ and $S^{\pm}$ has
  one of the forms illustrated in Figure~\ref{zscontainoney}

  \begin{figure}[h]
    \psfrag{e}{$e$}
    \psfrag{g}{$g$}
    \psfrag{fbar}{$\overline{f}$}
    \psfrag{R1}{$R_1$}
    \psfrag{Rn}{$R_n$}
    \psfrag{dots}{$\dotsb$}
    \psfrag{S0}{$S_0$}
    \centerline{%
      \includegraphics{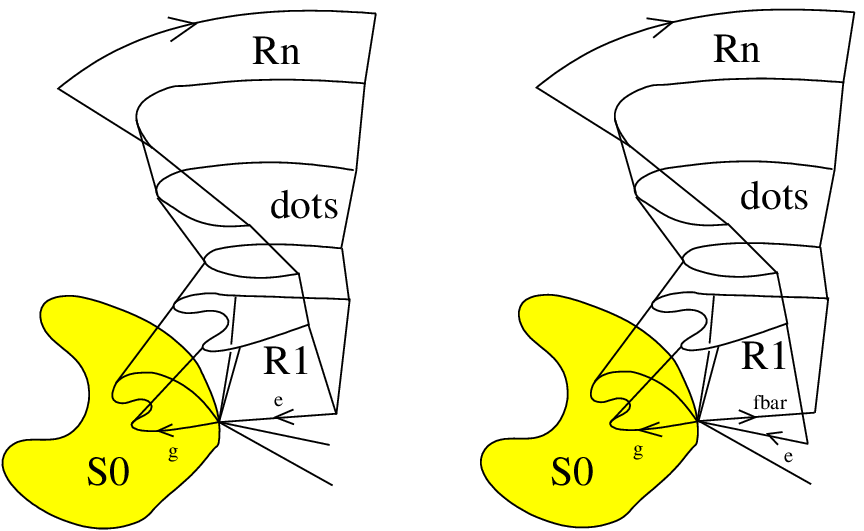}
    }
    \caption{}
    \label{zscontainoney}
  \end{figure}

  Orient each $\partial^{\pm}R_i$ so that $\partial^+R_i$ and
  $\partial^-R_{i-1}$ are oriented coherently, and so that under the
  map $Z\to N$ (recall that $Z$ is a connected component of a preimage
  of $N\in\mathcal{N}(X_H)$.) each $\partial^{\pm}R_i$ maps to the
  image of $e_{k(i)}g_i\overline{e_{k(i)}},$ in the first case, and
  maps to $e_{k(i)}g_i\overline{f_{l(i)}}$ in the second.

  In the second case, without loss, suppose that $\partial^+R$ and
  $\partial^+R_n$ are oriented the same way. By the same argument
  (Lemma~\ref{structurefortrees}) used to show that in a
  union of trees $Z,$ the graph of spaces associated to a maximal tree
  in $G(Z)$ is treelike, we know that $Z'_i$ is treelike, and in order
  for $Z$ to be treelike, we must have that $\partial^-R$ is attached
  to $S_0$ along a path $e_{k(i)}g'\overline{f_{l(i)}},$ otherwise the
  number of boundary components must strictly decrease. In $N,$
  $e_{k(i)}g\overline{f_{l(i)}}$ is homotopic to
  $e_{k(i)}g'\overline{f_{l(i)}},$ and since $Y\to N$ is $\pi_1$
  injective, we have that $g_i$ and $g'$ represent the same element of
  $\pi_1(Y)$. By Convention~\ref{chooserepresentatives} $g_i=g'$.

  In the first case, assume again that $\partial^+R_n$ and
  $\partial^+R$ are oriented the same way. Then for the same reason
  that $\partial^-R$ is attached to a path
  $e_{k(i)}g'\overline{f_{l(i)}}$ in the previous case, $\partial^-R$
  is attached along a path $e_{k(i)}g'\overline{e_{k(i)}}$ in
  $S_0$. As before, $g_i=g'$ because $\pi_1(Y)$ embeds in $H$.

  In either case, every rectangle in $G_i$ is attached along the same
  path $e_{k(i)}g_i\overline{e_{k(i)}}$ or
  $e_{k(i)}g_i\overline{f_{l(i)}}$. For each $i,$ introduce a relative
  graph $S_i\cong I$ with edge path $S_i\to S_0$ agreeing with
  $e_{k(i)}g_i\overline{e_{k(i)}}$ or
  $e_{k(i)}g_i\overline{f_{l(i)}},$ as the case may be. Each rectangle
  $R\in G_i$ which meets $S_0$ has attaching map $\partial^-R\to S_0$
  which factors through a map $\psi_R\colon\partial^-R\to S_i$. Let
  $Z_i$ be the graph of spaces $(G_i\setminus\set{S_0})\cup\set{S_i}$ with
  attaching maps $\psi_R$ for appropriate $R$ or attaching maps
  agreeing with the original attaching maps if $R$ does not meet $S_0$.

  We recover $Z$ by attaching each $Z_i$ along $S_i$ to $S_0$ via
  $h_i=e_{k(i)}g_i\overline{f_{l(i)}}$ (where $f_{l(i)}=e_{k(i)}$ in
  the first case).  Each $Z_i$ is a union of trees
  (Definition~\ref{unionoftrees}) which, if one were not treelike,
  would imply that $Z$ is not treelike.\mnote{by kappa stuff! write
    this down!}
\end{proof}

\par Before we proceed, perform a homotopy of $\widetilde{\varphi}$
which will simplify the analysis a little: Fix $Z$ containing a single
$Y,$ and let $\set{Z_i},$ $\set{S_i},$ and $S_0$ be as in the previous
lemma. The map $\widetilde{\varphi}$ carries $(Z,\partial Z)$ to some
$(N,\partial N)\in\mathcal{N}(X_H)$. Consider the map
$Z_i\to\Gamma_{Z_i}$ constructed after
Definition~\ref{unionoftrees}. Since each element of
$\YY(X_H)$ was chosen to be aspherical, by
Lemma~\ref{extensionlemma}, we may choose a homotopy of
$\widetilde{\varphi},$ supported on $Z_i,$ fixed on $\partial Z_i$ and
$S_i,$ and so that the restrictions $\widetilde{\varphi}\vert_{Z_i}$
factor through maps $(\Gamma_{Z_i},\partial\Gamma_{Z_i})\to(N,\partial
N)$.

\par Given $Z$ containing some $Y,$ Build a rose $R_Y$ with basepoint
$b$ and a petal $p_i$ for each $Z_i$\footnote{It makes no difference
if the collection $\set{Z_i}$ is empty or not.} Define a map
$\psi_Y\colon R_Y\to Y$ such that the $i$'th petal maps to the path
$g_i$ in $\pi_1(Y,b)$ such that the edge path $h_i\colon S_i\to S_0$
that $Z_i$ is attached to $S_0$ along is precisely the path
$e_{k(i)}g_i\overline{e_{l(i)}},$ for the appropriate $k(i),$ $l(i)$.

\par Let $\Gamma_{Z}$ be the graph $R_Y\cup_{b=\tau(e_i)} e_i$. There is
an obvious map $\Gamma_Z\to S_0,$ mapping $\Gamma_Z\supset e_i\to
e_i\subset S_0$ and $R_Y\to Y$ via $\psi_Y$. Attach each $Z_i$ along
$S_i$ to $\Gamma_Z$ along the path
$\widetilde{h_i}=e_{l(i)}p_i\overline{e_{k(i)}}$ to build a graph of
spaces $V(Z),$ $\partial V(Z)=\cup\partial Z_i\subset
V(Z)\cong\partial Z$. Then $Z$ is recovered by attaching $V(Z)$ to
$S_0$ by the map $\Gamma_Z\to S_0$. 

Define \[\mathcal{S}(V(Z))=\set{\Gamma_Z}\cup\bigcup(\mathcal{S}(Z_i)\setminus\set{S_i})\] The
relationship between $\set{S_i},$ $\Gamma_Z,$ $\set{Z_i},$ $V(Z)$ and $Z$
is illustrated in the following triple of pushouts.
\begin{figure}[h]
\centerline{%
  \xymatrix{%
    \coprod S_i\ar[r]\ar[d]\ar[drr] & \Gamma_Z\ar[d]\ar[dr] & \\
    \coprod Z_i\ar[r]\ar[drr] & V(Z)\ar[dr] & S_0\ar[d] \\
    & & Z
  }
}
\end{figure}

For each $Z_i$ we have the quotient map ${\pi_{Z_i}}\colon
Z_i\to\Gamma_Z$. The attaching map $h_i\colon S_i\to S_0$ factors
through $\widetilde{h_i}\colon S_i\to \Gamma_Z$. Since
$S_i\in\mathcal{S}(Z_i),$ by Lemma~\ref{extensionlemma}, $S_i$ embeds
in $\Gamma_{Z_i}$ under $\pi_{Z_i},$ thus there is an induced map
\[\overline{\overline{h_i}}\colon\img_{\Gamma_{Z_i}}(S_i)\to\Gamma_Z\]
For this reason we call the image of $S_i$ in $\Gamma_Z$ $S_i$ as
well. Define a graph $\underlying(V(Z))$ as the
pushout in Figure~\ref{underlyingvpushout}.

\begin{figure}[h]
\centerline{%
  \xymatrix{%
    \coprod Z_i\ar[r] & \coprod \Gamma_{Z_i}\ar[dr] & \\
    \coprod S_i\ar[u]\ar[ur]\ar[dr] & & \underlying(V(Z)) \\
    & \Gamma_Z\ar[ur] &
  }
}
\caption{The diagram defining $\underlying(V(Z))$}
\label{underlyingvpushout}
\end{figure}

\par Since the restriction of $\widetilde{\varphi}$ to $Z_i$ factors
through $\Gamma_{Z_i}$ (after the homotopy provided by
Lemma~\ref{extensionlemma}), and $\underlying(V(Z))$ is the union of
$\Gamma_{Z_i}$ and $\Gamma_Z$ along $S_i,$ and since $Z$ is recovered
by attaching $V(Z)$ along $S_i$ to $S_0,$ the restriction of
$\widetilde{\varphi}$ to $V(Z)$ factors through the restriction to
$\underlying(V(Z))$. On the other hand, $S_0$ is simply
$Y_Z\cup\Gamma_Z,$ identified along $R_Z$. Hence,
$\widetilde{\varphi}\vert_Z$ is the composition of the projection of
$V(Z)$ to $\underlying(V(Z))$ followed by the restriction to $S_0$.

\par One key property of $V(Z)$ is that $\partial V(Z)$ is precisely
$\partial{Z}$. For $Z$ not containing $Y\in\YY(X_G),$ set
$V(Z)=Z$. Let $\mathcal{V}$ be the collection of all $V(Z)$.  Recall
the definition of $\mathcal{B}$. For each $B\in\mathcal{B}$ there are
two inclusions, each a homeomorphism with a boundary component of
$\ZZ_{\bullet}$. Let $\overline{B}$ be a point, as in the
construction of $Z'$. Since each $B$ maps to some boundary component
of some $Z,$ and each boundary component of a $Z$ is a boundary
component of $V(Z),$ and since each
$\widetilde{\varphi}\vert_{B|V(Z)}$ factors through
$\set{\overline{B}|\underlying(Z)},$ we may sensibly form the pushouts
in Figure~\ref{pushoutsforX}.

\begin{figure}[h]
  \centerline{%
    \xymatrix{%
      \coprod B^{(0)}\ar[r]\ar[d] & \coprod\mathcal{S}_{\bullet}(V(Z))\ar[r]\ar[d] & \Gamma(X)\ar[d] \\
      \coprod B\ar[r]\ar[d]       & \coprod V(Z)\ar[r]\ar[d]                       & X\ar[d] \\
      \coprod \overline{B}\ar[r]  & \coprod\underlying(V(Z))\ar[r]     & \underlying(X)
    }
  }
  \caption{Pushouts defining $\Gamma(X),$ $X,$ and
  $\underlying(X)$. Vertical arrows in the top row are inclusions,
  vertical arrows in the bottom row are projections.}
  \label{pushoutsforX}
\end{figure}

\par We see from the construction that $\Gamma(X)$ and
$\underlying(X)$ are graphs, and, since point preimages of
$V(Z)\to\Gamma_Z$ are connected, that
$\pi_1(X)\to\pi_1(\underlying(X))$ is onto. We can give an alternate
description of $X$ as
$\compcomp(X)\cup_{\set{\gamma_{0,i}}}\set{M_{0,i}}$. For each component
$Y\in\YY(X_G),$ there is exactly one $Z\in\ZZ$
containing it. Let $R_Y$ be the rose contained in $\Gamma_Z$. Then we
recover $\xtilde$ by gluing each $R_Y\subset\underlying(V(Z))$ to $Y$
via the map $\psi_Y$. The images of each $R_Y$ are disjoint under
$X\to\underlying(X)$. It follows easily from the fact that each
$\Gamma_{Z_i}$ embeds in $\Gamma_Z$ under $\pi_Z$ that the map
$\Gamma(X)\to\underlying(X)$ is an immersion.

\par Each $\underlying(V(Z))$ is a graph with distinguished boundary
vertices. Declare that every point of $\underlying(V(Z))$ whose
complimentary components number at least three to be vertices of
$\underlying(X),$ and let $V_1,\dotsc,V_{n_Z}$ be the leaves of
$\mathcal{F}(Z)$ (These are the only leaves which don't necessarily
have neighborhoods which are products.) which map to the vertices of
$\underlying(V(Z))$. Each such leaf is a finite graph. Setting all
such leaves, along with components of $\mathcal{B}_{\bullet}\subset
X,$ to be vertex spaces, it is easily seen that $X$ is a 2-covered
graph of spaces.

\begin{lemma}
  $\chi(\Gamma(X))=\chi(\underlying(X))$
\end{lemma}

\begin{proof}
If $\Gamma$ is a graph with distinguished valence one vertices
$\partial\Gamma,$ define
\[\kappa(\Gamma,\partial\Gamma)=-\chi(\Gamma)+\frac{1}{2}\#\partial\Gamma\] 
as was previously done for graphs with no loops. Observe that
\[\sum_{Z\in\ZZ}\sum_{S\in\mathcal{S}(V(Z))}\kappa(S,\partial
S)=-\chi(\Gamma(X))\] and
\[\sum_{Z\in\ZZ}\kappa(\underlying(V(Z)),\partial\underlying(V(Z))=-\chi(\underlying(X))\]
Thus we only need to check that $\kappa(\underlying(V(Z)),\partial
V(Z))=\sum_{S\in\mathcal{S}(V(Z))}\kappa(S,\partial S)$ for all
$Z$. This follows easily from
$\chi(\underlying(V(Z)))=\chi(\Gamma_Z),$ $\Delta_q^-(Z_i)=0,$ and the
fact that $\partial S_i$ maps to distinct boundary components of
$\Gamma_{Z_i}$ as in the computation carried out in
Lemma~\ref{deltapforz3}.
\end{proof}




\section{Moves on graphs of spaces}%
\label{movesongraphsofspaces}

\par A $2$-covered graph of spaces is generally an ugly beast, but can
be convert ed to a more amenable object by \term{folding},
\term{reducing}, and \term{collapsing}. We handle them in reverse
order.

\begin{definition}[Collapse] 
  If $X$ is a graph of spaces and $e$ is an edge of
  $\underlying(X)$ with $\tau(e)\neq\iota(e),$ and if
  $\tau\colon E\to V_{\tau(e)}$ is an embedding, then we can
  collapse $X$ to $X_{E}$ by crushing the edge space
  $e\times E$ to $\iota(E)$. In the topological realization of
  $X,$ collapse $E\times\mathrm{I}$ to
  $E\times{\set{0}}$. The resulting vertex is
  $E_{\iota(e)}\cup E_{\tau(e)}/\tau(w)\sim\iota(w),$
  $w\in E$. The edge maps incident to the new vertex are still
  immersions, and it's easy to check that the quotient map is a
  homotopy equivalence which respects
  $\pi_1(\Gamma(X))\to\pi_1(X)$. \mnote{Mark is right. This is
    confusing.}
\end{definition}

\mnote{the collapses here aren't the full list of possible collapses}

\begin{definition}[Weight]
  The weight of a graph is the number of edges.
\end{definition}

\begin{definition}[Reduced]
Some graphs of spaces admit trivial simplifications. For instance, the
topological realization of a graph with a valence two vertex can be
given a simpler description by un-subdividing an edge. A similar
statement holds for our 2-covered graphs: If $V$ is a vertex space
in $X$ and $E_1$ and $E_2$ are the only incident edges,
then if both maps $E_i\to V$ are graph isomorphisms, then
$X$ is \term{reducible}. By collapsing one of the incident edges,
the number of reducible vertices strictly decreases. If $X$ has no
reducible vertices, and all valence one vertices have nonzero weight,
then it is \term{reduced}.
\end{definition}

\begin{definition}[Folding]
  Fix a $2$-covered graph of spaces $X$. Given a set of edges,
  indexed by $J,$ we define a new graph of spaces $X_J,$ called
  a \term{fold} of $X$. We say that $X_J$ is obtained from
  $X$ by \term{folding}.

  How to fold: Let $V$ be a vertex of a graph of spaces $X$ as
  above. Let $\set{(E_i,\tau_i)}_{i\in I}$ be the oriented edges
  whose terminal vertex is $V$. For $J\subset I,$
  define \[V_J=\bigcup_{j\in J}\tau_j(E_j)\]

  \par Let $\set{V_{J,p}}_{p=1..l_p}$ be the connected components of
  $V_J$ and $\set{V_{I\setminus J,q}}_{q=1..l_q}$ the connected
  components of $V_{I\setminus J},$ and $\set{E_{J,r}}_{r=1..l_r}$
  the connected components of $V_J\cap V_{I\setminus J}$.

  \par For each index $p,$ $q,$ $r,$ introduce new vertices $v_{J,p},$
  $v_{I\setminus J,q},$ and oriented edges $e_{J,r}$ with
  $\tau(e_{J,r})$ the member of $\set{V_{J,p}}$ that
  $E_{J,r}$ is contained in, and $\iota(e_{J,r})$ the member of
  $\set{V_{I\setminus J,q}}$ that $E_r$ is contained in.
  Define $\iota\colon\overline{e_{J,r}}\to\cdots$ to be
  $\tau:e_{J,r}\to\cdots,$ where $\cdots$ represents the appropriate
  component $V_{J,p}$ or $V_{I\setminus J,q}$.

  \par This data, along with the (undisturbed) data from the rest of the
  graph of spaces $X$ defines a new graph of spaces (in the
  2-Covered sense) $X_J$ with the vertex space $V$ split apart.
\end{definition}

\begin{figure}[tb]
\psfrag{AcupB}{{\tiny $E_1\cup E_2$}}
\psfrag{A}{{\tiny $E_1$}}
\psfrag{B}{{\tiny $E_2$}}
\psfrag{AdeltaB}{{\tiny $E_1\triangle E_2$}}
\centerline{%
  \includegraphics[scale=0.6]{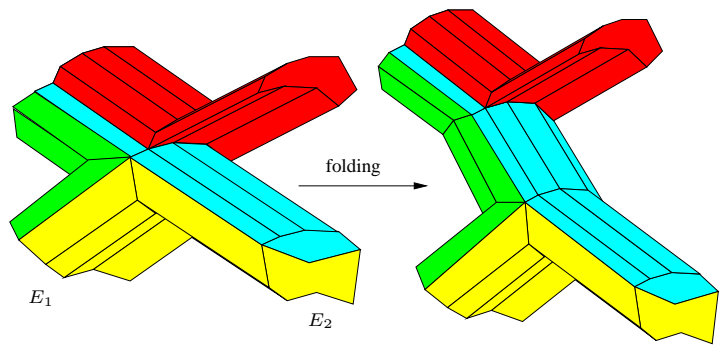}
}
\caption{Folding edges $E_1$ and $E_2$ together to simplify $X$.}
\label{splittingfigure}
\end{figure}

Folding is illustrated in Figure~\ref{splittingfigure}. Note that $J$
may consist of a single element, yet the split space may still be
distinct from the original space. Also, beware that it's possible for
the underlying graph's complexity to increase: the subgraph of
$\underlying(X_J)$ spanned by $e_{J,r}$ may not be a tree.

\begin{definition}[Unfoldable]
  A vertex $v\in\underlying(X)$ is \term{unfoldable} if for all $J\subset
  I,$ where $I$ is the indexing set of the incident edges $E_i,$
  one of
  \[V_J\simeq\coprod_{j\in J} E_j\quad\mbox{or}\quad V_{I\setminus J}\simeq\coprod_{i\in I\setminus J} E_i\]
  holds. If a vertex isn't unfoldable, then it is \term{foldable}.
\end{definition}

\par Unfoldable vertices are particularly nice. Not only do they fall
into two basic simple types, folding an unfoldable vertex doesn't
change the graph of spaces.

\begin{lemma}[Structure of Unfoldable Vertices]
\label{unfoldablestructure}
A reduced unfoldable vertex $v$ has the form
\begin{itemize}
  \item There is a distinguished edge $e_0$ adjacent to $v$. The rest
    of the edges $e_1,\dotsc,e_m$ are undistinguished.
  \item The map $E_0\immerses V$ is not an embedding.
  \item The maps $E_i\immerses V,$ $i\neq0,$ are embeddings
    with pairwise disjoint images.
\mnote{Removed last entry. Does that make it more confusing? Check if
  some other proof relies on it.}
\end{itemize}
or 
\begin{itemize}
\item $v$ has valence three and all incident edge maps are
  embeddings. The image of every incident edge space meets every
  other. There is a vertex $w$ of $V$ which is in the image of
  every incident edge space.

A fold of an unfoldable vertex in $X$ recovers $X$.
\end{itemize}
\end{lemma}

\begin{definition}
An unfoldable vertex with a distinguished edge $e_0$ such that
$E_0\immerses V$ is not an embedding is \term{degenerate}. A
vertex that has valence three and whose incident edges embed is
\term{nondegenerate}.
\end{definition}

\begin{proof}[Proof of Lemma~\ref{unfoldablestructure}]
\par Let $v$ be an unfoldable vertex.
\par If an incident edge $E_0\immerses V$ isn't an embedding, then 
it's clear we're in the first case of the lemma. Take $J=\set{0}$. Then
the graph covered by the remaining incident edge graphs is
homeomorphic to their disjoint union.

\par Thus we need to show that the second case of the lemma holds,
assuming every incident edge map is an embedding. Suppose that the
valence of $v$ is at least four. Either there is a chain of incident
edge graphs $E_i,$ $i=1,2,3,4,$ such that
$\img(E_i)\cap\img(E_{i+1})\neq\emptyset$ or there is an
incident edge $E_1$ whose image meets every other incident edge
graph. In the first case, we may take $J=\set{1,2}$.

\par In the second case, If $E_i,$ $i,j,k\neq 1,$ $i,j,k$
distinct, whose images meet $E_1,$ then they must have disjoint
images since there is no chain of length four. For example, if
$\img(E_2)\cap\img(E_3),$ then the sequence
$(E_4,E_1,E_2,E_3)$ is a chain of length four. Since
$V$ is connected, there is an edge $f$ of $V,$ contained in
the image of $E_1,$ which isn't covered by any $E_i,$ $i\neq
1,$ thus $f$ is covered twice by $E_1,$ a contradiction.

Let $V'=\img(E_1)\cup\img(E_2)$. If $f$ is an edge of
$V$ meeting $V'$ and $f$ isn't contained in $V'$ then $f$
is covered by $E_3$. The endpoint of $f$ contained in $V'$ is
contained in the image of every incident edge space.
\end{proof}

\par If $E_{J,r}$ is the set of edges introduced by folding a set of
incident edges $\set{E_j}_{j\in J},$ then the original graph of spaces
is recovered by collapsing $\set{E_{J,r}}$.

\par Let $\varphi\colon X_J\to X$ be the collapsing map. Then
$\varphi$ is a homotopy equivalence. Let $\Gamma$ be a connected
component of $\Gamma(X)$.  If $\Gamma_J$ is the associated
connected component of $\Gamma(X_J),$ $\varphi_J$ is the
collapsing map restricted to $\Gamma_J,$ and $\varphi_{U}$ is induced
map on underlying graphs, then

\centerline{%
  \xymatrix{%
    \Gamma_J\ar[r]^{\imath}\ar[d]^{\varphi_J} &  X_J\ar[r]\ar[d]^{\varphi} & \underlying(X_J)\ar[d]^{{\varphi_{U}}} \\
    \Gamma\ar[r]^{\imath}   &  X\ar[r]   & \underlying(X)
  }
}

\noindent commutes. Collapsing restricted to the horizontal subgraph
crushes forests, thus \mnote{if this is okay, then i need to define
  collapsing more precisely} $\varphi_J$ is a homotopy equivalence,
${\varphi_U}_*$ is an epimorphism, and the unlabeled arrows are the
natural epimorphisms \[\mbox{graph of spaces}\to\mbox{underlying
  graph}\]

\par Given a graph of 2-covered graphs, there is a reduced space
$X\mapsto X^R$ obtained by trimming trees and removing all
valence two vertices for which both incident edge maps are graph
isomorphisms.

\section{Simplifying graphs of spaces}
\label{simplifying}
Under certain favorable conditions a folded space admits further
simplification. There is a complexity, which, when minimized through
folding and collapsing, gives an optimal graph of spaces equivalent to
a given one. The structure of the vertex and incident edge spaces of a
space minimal with respect to this complexity is considerably simpler
than that of a nonminimal graph of spaces.

\begin{definition}[Complexity of Graphs of Spaces]
\label{complexityofgraphsofspaces}
  Let $k(X)$ be the maximal valence of a vertex in
  $\underlying(X),$ $\mathrm{m}_l(x)$ the number of vertices of
  valence $l,$ $\mathrm{m}_2^{red}$ is the number of reducible valence two
  vertices, and $\mathrm{m}_2^{deg}(X)$ is the number of
  degenerate valence two vertices in $\underlying(X)$. If $X$
  is reduced, and $k(X)\geq3,$ then the complexity of $X$ is
  the tuple \[c(X)= (-\betti(\underlying(X)),
  k(X),m_{k(X)}(X), \dotsc,m_3(X),
  m_2^{red}(X),-m^{deg}_2(X))\] If $k(X)=2,$ then the entries
  $\mathrm{m}_{k(X)}(X),\dotsc,m_3(X)$ don't appear. The
  order is the lexicographic one.
\end{definition}

\par $X'$ is obtained from $X$ by folding if there is a
sequence of folds $X= X_0\to X_1\to\dotsb\to X_k= X'$.
Let $\mathrm{Folds}(X)$ be the set of graphs of spaces which can
be obtained by folding.


\begin{lemma}[Minima of $c$]
\label{foldtovalencethree}
Let $X\in\mathrm{Folds}(Y)$. If $c(X)$ is minimal then
$X$ is reduced and every vertex of $X$ is unfoldable.
\end{lemma}

\par
\begin{figure}[tb]
\centerline{%
  \includegraphics{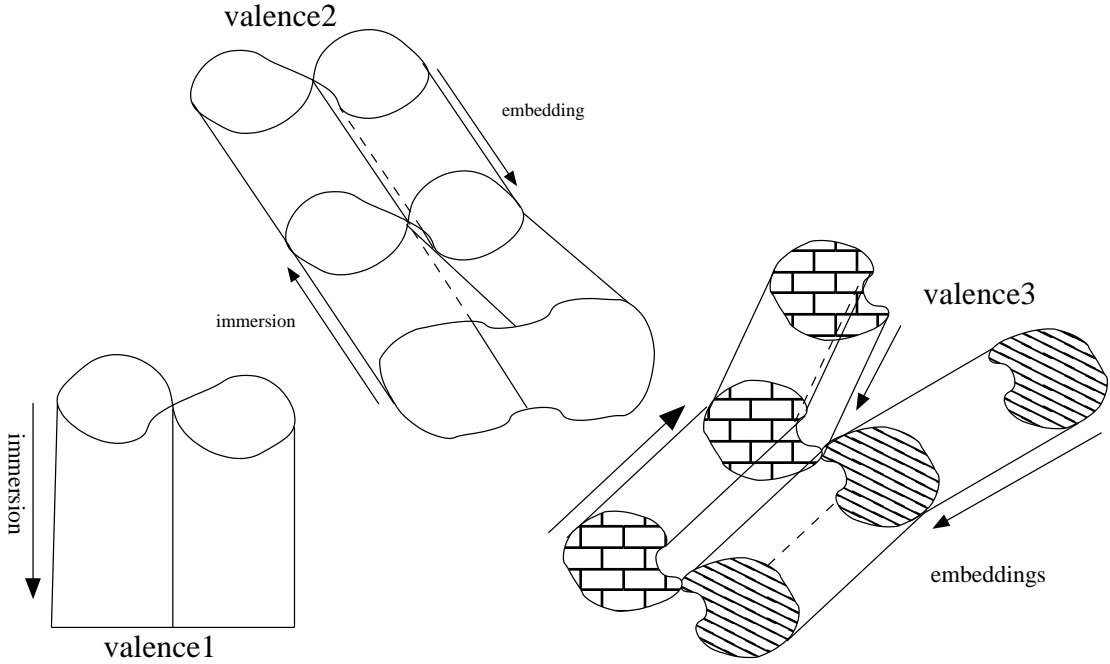}
}
\caption{Illustration of Lemma~\ref{foldtovalencethree}.}
\label{valencethreeillustration}
\end{figure}

 The conclusions of Lemma~\ref{foldtovalencethree} are illustrated
in Figure~\ref{valencethreeillustration}. 

\begin{proof}
If $X$ can be reduced, then reducing decreases $c$. The idea here
is that folding takes a vertex in $\underlying(X)$ and blows it up
to a bipartite graph. Then either $\betti(\underlying(X))$
increases or the graph is a tree. If it's a tree, then either $k$ or
$\mathrm{m}_k$ must decrease unless the vertex is unfoldable. Also
note that $c$ takes only finitely many values on
$\mathrm{Folds}(Y)$..

If $\underlying(X)$ has a foldable vertex then $X$ isn't a
minimum of $c$: Let $v$ be a foldable vertex, and let
$\set{E_j}_{j\in J}$ such that neither
\[V_J\simeq\coprod_{j\in J} E_j\quad \mbox{nor}
\quad V_{I\setminus J}\simeq\coprod_{i\in I\setminus J} E_i\eqno{(\diamondsuit)}\]
holds.

\par Suppose $J=\set{1,2}$. Let $v_J$ be the 
additional vertex corresponding to $V_J$. Let $v_1,\dotsc,v_q$ be
the vertices corresponding to connected components of
$V_{I\setminus J},$ and let $e_1,\dotsc,e_r$ be the edges
corresponding to connected components of $V_J\cap V_{I\setminus
J}$. Let $s$ be the valence of $v$. 

If $r>q$ then
$\betti(\underlying(X_J))>\betti(\underlying(X)),$ thus
$c(X)$ isn't minimal.

Thus we may assume that $r=q$. There are $s-2$ edges incident to
$v_1,\dotsc,v_q$. Each vertex $v_j,$ $j=1..q$ has valence at most
$s-1$ and $v_J$ has valence at most $s$. If $v_J$ has valence $s,$
then $r=q=s-2$ and the edges $E_i,$ $i\in I\setminus J,$ have
pairwise disjoint images, therefore at least one must be an immersion
but not an embedding, increasing $\mathrm{m}_2^{deg}$. 

\par If $J=\set{1},$ a similar argument works. If the valence of $v_J$ is $s,$ 
then $r=q=s-1,$ implying that the images $\img(E_i),$ $i\in
I\setminus J,$ are pairwise disjoint, therefore at least one of them
is an immersion but not an embedding, increasing
$\mathrm{m}_2^{deg}$. If a vertex $v_i$ has valence $s,$ then $r=q=1,$
$v_J$ has valence two, and $\mathrm{m}_2^{deg}$ increases. Thus, in all
cases, if $X$ has a foldable vertex, then $c$ can be decreased by
folding.

\par If $J$ has more than two elements, then there is a subset of $J$ with 
either one or two elements which satisfies $(\diamondsuit)$.
\end{proof}


\begin{lemma}[Euler Characteristic Lemma]
  \label{chiinequality} Let $X$ be a 2-covered graph of
  spaces. Then
  \[\chi(\Gamma(X))\leq\chi(\underlying(X))\]
\end{lemma}

%

\begin{proof}[Proof of Lemma~\ref{chiinequality}]
  \par First, assume that $X$ is a minimum of $c,$ i.e.,
  $\mathrm{Folds}(X)=\set{X}$. We handle the different valence
  vertices of $\underlying(X)$ on a case-by-case basis. The Euler
  characteristic of a graph can be computed by adding the
  ``curvatures'' of its vertices:

 \[ \chi(\Gamma)=\sum_{v\in\Gamma^{(0)}} \kappa(v),\quad
 \kappa(v)=1-\frac{1}{2}\mathrm{valence}(v)\]

 For each vertex $v$ of $\underlying(X),$ let $\widetilde{v}$ be
 the set of vertices of $\Gamma(X)$ which map to $v$. Let
 $\kappa({\widetilde{v}})=\sum_{w\in\widetilde{v}}\kappa(w),$ and so
 $\chi(\Gamma(X))=\sum \kappa({\widetilde{v}})$.

 There are two cases to consider. Recall the structure of unfoldable
 vertices from Lemma~\ref{unfoldablestructure}.

 {\bf $\boldsymbol{v}$ is degenerate:} Let $k+1$ be the valence of
 $v$. Let $E_0$ be the immersed edge graph, and $E_i,$
 $i=1..k,$ the embedded edge graphs. Every vertex of $V$ is the
 image of at least two distinct vertices of $\cup E_i^{\circ},$
 hence has valence at least two.
 
 \par Let $V_1$ be the union of edges of $V$ covered twice by
 $E_0$. The vertex graph is the union
 $\cup_{i\neq0}\img(E_i)\cup V_1$.
 
 \par Suppose $V_1$ is nonempty. Let
 $V_1^{\circ}=(\img(E_0)\setminus\cup_{i\neq
 0}\img(E_i))\subset V_1$ be the subgraph of $V$ covered by
 $E_0$ but \emph{not} covered by any other edge space. Then
 $V_1^{\circ}$ is the interior of $V_1$. \mnote{FIX: this
 definition isn't great, but it's accurate} A vertex in
 $V_1^{\circ}$ contributes a vertex with valence at least two to
 $\Gamma(X)$: If $w$ is a vertex in $V_1^{\circ},$ it has an
 incident edge $f,$ $\tau(f)=w,$ which, since $\tau$ is an immersion,
 is the image of two \emph{distinct} oriented edges $f_1$ and $f_2,$
 $\tau(f_1)\neq\tau(f_2),$ from $E_0$. The terminal vertices of
 $f_i$ map to $w$ under $\tau,$ thus the vertex $w,$ regarded as a
 vertex of $\Gamma(X),$ has valence at least two in
 $\Gamma(X)$.
 
 \par We now handle the vertices $(\cup_{i\neq
 0}\img(E_i))\cap V_1$.  Edges not contained in $V_1$ are
 each covered once by $E_0$ and once by $\cup_{i\neq0} E_i$.
 Since $V$ is connected, there are oriented edges
 $f_i\subset V_1$ which meet $\img(E_i)$ at their terminal
 vertices $w_i$. Each $f_i$ is the image of distinct oriented edges
 $f^i_1\neq f^i_2\subset E_0$ with distinct (since $E_0$
 immerses in $V$) terminal vertices $w^i_1$ and $w^i_2$. Thus
 $w_i$ is the image of a vertex in $E_i,$ and the image of two
 vertices in $E_0,$ hence $w_i$ has valence at least three vertex
 in $\Gamma(X)$.
 
 \par If $V_1$ is empty, then there is a vertex
 $w\in V=\img(E_0)=\img(E_1)$ which is the image of two
 vertices in $E_0$. Then $w$ is the image of a vertex in $E_1$
 as well. Then $w,$ regarded as a vertex of $\Gamma(X),$ has
 valence at least three.
 
 \par In all cases $\kappa({\widetilde{v}})\leq -1/2\cdot
 k<\kappa(v)=(2-(k+1))/2=1/2-1/2\cdot k$. The inequality is strict.
 
 {\bf $\boldsymbol{v}$ is nondegenerate:} There are three edges
 incident to $V$: $E_1,$ $E_2,$ and $E_3,$ and all
 incident edge maps are embeddings.
 
 \par Suppose $E_1$ is a point. Since
 $\tau_2\colon E_2\to V$ and $\tau_3\colon E_3\to V$ are
 both embeddings, and every edge is covered once by $E_2$ and once
 by $E_3,$ both are surjective. A surjective immersion of graphs
 is an isomorphism, hence\mnote{clumsy hence both$\dotsc$} both maps are
 graph isomorphisms. Let $w$ be the image of $E_1$. The incident
 edge maps map $w_2\in E_2$ and $w_3\in E_3$ in $E_2$ to
 $w$. Since $w$ is also the image of $E_1,$ it is a valence three
 vertex of $\Gamma(X)$. Every other vertex in $\widetilde{v}$ has
 valence two in $\Gamma(X),$ hence contributes nothing to
 $\kappa({\widetilde{v}})$. Thus
 $\kappa({\widetilde{v}})=\kappa(v)=-1/2$.
 
 \par We're left with the possibility of three nontrivial embeddings,
 i.e., $E_{1\vert2\vert3}$ aren't points. Every vertex of
 $\widetilde{v}$ has, by the previous arguments, valence at least
 two. Since the incident edge maps are embeddings, every vertex is
 covered at most once by each incident edge, i.e., every vertex in
 $\widetilde{v},$ regarded as a vertex of $\Gamma(X),$ has valence
 at most three. Since $v$ is nondegenerate there exists a point of
 triple intersection, hence
 
 \begin{itemize} 
 \item $\kappa({\widetilde{v}})\leq\kappa(v)=-1/2$
 \item If $\kappa({\widetilde{v}})=\kappa(v)$ then there is exactly one
   point of triple intersection of incident edge graphs, otherwise
   $\kappa({\widetilde{v}})\leq -1<\kappa(v)=-1/2$
 \end{itemize}
 
 \par By the cases above, we conclude the inequality 
 \[\chi(\Gamma(X))=
 \sum \kappa({\widetilde{v}}) \leq \sum \kappa(v) =
 \chi(\underlying(X))\eqno{(\spadesuit)}\] for minima of $c$.
 
 \par If $X$ isn't a minimum of $c,$ then let $X_c$ be a
 member of $\mathrm{Folds}(X)$ with minimal complexity. Since
 $\betti(\underlying(X_c))\geq\betti(\underlying(X))$
 \[\chi(\underlying(X))\geq\chi(\underlying(X_c))\geq\chi(\Gamma(X_c))=\chi(\Gamma(X))\eqno{(\clubsuit)}\]
\end{proof}



\par We'll be interested in graphs of spaces whose horizontal
subgraphs have the same Euler characteristic as their underlying
graphs. When this happens, the space can be folded so that all
vertex spaces have the simplest form possible.

\begin{lemma}[$\chi(\Gamma(X))=\chi(\underlying(X))$]
\label{sameeulercharacteristic}
Suppose $X$ is a minimum of $c$ and
$\chi(\Gamma(X))=\chi(\underlying(X))$. Then every vertex has
valence three. If $V$ is a vertex with incident edge spaces
$E_i,$ $i=1,2,3,$ then $\cap E_i$ is a single point.
\end{lemma}


\begin{proof}
\par Suppose $\chi(\Gamma(X))=\chi(\underlying(X))$.  By
$(\clubsuit),$ every minimum $X_c$ of $c$ obtained by folding
satisfies
$\betti(\underlying(X_c))\geq\betti(\underlying(X))$. If this
inequality is strict, then, by Lemma~\ref{chiinequality},
$\chi(\Gamma(X))<\chi(\underlying(X))$. Thus
$\chi(\Gamma(X))=\chi(\underlying(X)),$ and for every minimum
$X_c,$ $\chi(\Gamma(X_c))=\chi(\underlying(X_c))$.

\par Let $X_c$ have minimal $c$ out of all members of
$\mathrm{Folds}(X)$. If $X_c$ had a degenerate vertex, then
the inequality $(\spadesuit)$ would be strict, thus every vertex is
unfoldable, has valence three, and is nondegenerate.  By the argument
used to prove Lemma~\ref{chiinequality}, there is exactly one point of
triple intersection of edges incident to every vertex graph.
\end{proof}

\mnote{this underwent small reconstruction. Is it still okay?}

\begin{lemma}
  \label{tripleintersection}
  A graph $V,$ 2-covered by connected subgraphs
  $E_i\into V,$ $i=1,2,3,$ such that $\cap\img(E_i)$ is a
  single vertex $w,$ has one of the following forms:
  \begin{itemize}
    \item $E_i$ are all points.
    \item $E_1$ is a point, and
      $E_{2\vert3}\cong V$. $\mathscr{W}(E_{2\vert3})>0$.
    \item $V=\img(E_1)\vee_w\img(E_2),$
      $E_3\cong V$. $\mathscr{W}(E_{1\vert2\vert 3})>0$.
    \item $V$ is the union of three subgraphs $V_{1\vert 2\vert 3}$ which
      meet at a single vertex
      $w\in V\cong\bigvee_w V_i$. $E_i\cong V_{i+1}\vee_w V_{i+2}$. $\mathscr{W}(E_{1\vert 2\vert 3})>0$.
\end{itemize}
\end{lemma}

\begin{definition}[Separable, Trivial, Splittable]
A vertex $v$ of a graph of spaces such that $V$ satisfies
Lemma~\ref{tripleintersection} is called \term{separable}. If $V$
satisfies one of the first two bullets $v$ is
\term{trivial}. Otherwise it is \term{nontrivial}. If $V$
satisfies the third bullet $v$ is \term{splittable}. If $V$
satisfies the fourth, $v$ is separable, but \term{unsplittable}.
\end{definition}

\begin{figure}[b]
\psfrag{AcapB}{{\tiny $E_1\cap E_2$}}
\psfrag{BcapC}{{\tiny $E_2\cap E_3$}}
\psfrag{AcapC}{{\tiny $E_1\cap E_2$}}
\psfrag{p}{{\tiny $p$}}
\psfrag{ell}{{\tiny $\dotso$}}
\centerline{%
  \includegraphics[scale=0.5]{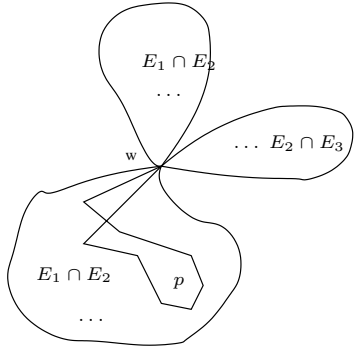}
}
\caption{$w$ separates $V$. $p$ is a path from the proof of Lemma~\ref{tripleintersection}.}
\label{valencethreeseparates}
\end{figure}

\begin{proof}[Proof of Lemma~\ref{tripleintersection}]
\par Let $w=w_v$ be the sole point of triple intersection. Let
$\mathcal{P}$ be the set of edge-paths starting at $w$ and that
terminate if they meet $w$ again. Let $\set{E_i}_{i=0,1,2}$ be the
edge graphs incident to $V$. We divide $\mathcal{P}$ into three
subclasses $\mathcal{P}_{j,k},$ $j\neq k$. A path $p$ lies in
$\mathcal{P}_{j,k}$ if the image of $p$ is contained in
$\img(E_j)\cap\img(E_k)$. Let
$V_i=\bigcup_{p\in\mathcal{P}_{i+1,i+2}}\img(p)$. At most one of
$V_i$ can be a point. Every point/edge of $V$ lies in one of
$V_i,$ which all meet at $w,$ the sole point of triple
intersection. An incident edge $E_i$ is then isomorphic to
$V_{i+1}\cup V_{i+2}$.  
\end{proof}

\par See Figure~\ref{valencethreeseparates} for an illustration. In
virtue of Lemmas~\ref{sameeulercharacteristic}
and~\ref{tripleintersection} we make the following definition.

\begin{definition}
\label{separable}
  If $X$ satisfies $\chi(\underlying(X))=\chi(\Gamma(X))$
  and $X$ is a minimum of $c,$ since every vertex of
  $\underlying(X)$ is separable, we say that $X$ is
  \term{separable}.
\end{definition}

\section{Separable graphs of spaces}
\label{separablespaces}
\par In this section we consider only separable graphs of spaces.

\par The next two lemmas give us the means to analyze the separable
graphs of spaces. 

\begin{definition}
  \label{irreducible}
  A graph of spaces is \term{irreducible} if it has no trivial edge
  spaces, i.e., there are no ``obvious'' free product decompositions
  of its fundamental group. The removal of interiors of weight $0$
  edges and leftover vertices from $X$ yields graphs of spaces
  $X_i$ which are the \term{irreducible components} of $X$.
\end{definition}


\par Every vertex of a separable graph of spaces turns each edge space
into a (possibly trivial) wedge of subgraphs. We would like to push
this structure around the graph of spaces to give each edge graph the
coarsest treelike structure compatible with all decompositions forced
upon it. This lets us express a two-covered graph of spaces as a union
of \term{cylinders}. Under certain circumstances a graph of spaces can
be repeatedly collapsed and folded to what is essentially a wedge of
cylinders. 

\par We start by defining the cylinders of a graph of spaces
$X$. Roughly speaking, a cylinder $C$ is a graph of spaces
whose underlying graph is a circle, has a map to $X$ compatible
with edge maps, and if the map $C\to X$ factors through a
similar such map $C'\to X,$ then $C\cong C'$.

\begin{definition}[Cylinder]
\label{cylinder}
  A graph of spaces is a \term{cylinder} if its underlying graph is a
  circle and has only reducible valence two vertices. A cylinder is
  homeomorphic to the mapping torus of a combinatorial automorphism of
  a graph.

  Let $\mathbb{S}(X)$ be a set of indivisible (not factoring
  through a covering map $ S^1\to S^1$), unoriented, closed, immersed
  edge paths in $\Gamma(X)$ uniquely representing every conjugacy
  class of maximal cyclic subgroup of $\pi_1(\Gamma(X))$ as an
  immersion $\iota\colon S^1\to\Gamma(X)$. There is an immersion
  $\mathbb{S}(\iota)\colon\mathbb{S}(X)\immerses\Gamma(X)$. A
  graph of spaces $X$ is a union of annuli and \mobius\ bands
  $\set{A_j}$ and $\Gamma(X)$ along boundary maps
  $\varphi_j\colon\partial A_j\immerses\Gamma(X)$. Each annulus
  is a union of squares and the map $\varphi_j$ is a pair (or a
  singleton, if $A_j$ is a \mobius\ band) of edge paths in
  $\Gamma(X)$. The maps $\varphi_j$ factor through
  $\mathbb{S}(X),$ i.e, there are lifts
  \[\widetilde{\varphi_j}\colon\partial A_j\immerses\mathbb{S}(X)\]
  such that
  $\mathbb{S}(\iota)\circ\widetilde{\varphi_j}=\varphi_j$. This is
  because all edge maps $E\to V$ are immersions, hence the maps
  $\partial A_j\to\Gamma(X)$ are immersions.

  A graph of spaces $X$ is the union
  $\Gamma(X)\cup_{\varphi_j}A_j$. Define a new graph of spaces
  $\overline{X}$ to be
  $\mathbb{S}(X)\cup_{\widetilde{\varphi_j}}A_j$. The set of cylinders
  of $X,$ denoted $\Cyl(X),$ is the collection of connected components
  of $\overline{X}$ containing an annulus or \mobius\ band.
  
  The boundary of a cylinder $C,$ $\partial_{X}C,$ is the subgraph
  of $\Gamma(C)$ corresponding to elements of $\mathbb{S}(X)$
  whose images are contained in $\compcomp(X)$. The boundary map
  $\mathbb{S}(\iota)\vert_{\partial_{X}C}$ is denoted
  $\varphi_{C}$. The inclusion map of a cylinder (which isn't
  really an inclusion, but we ignore this technicality)
  $C\to X$ is denoted $\psi_{C}$.

  The space $X$ is recovered by forming the quotient space
  $\compcomp(X)\cup_{\varphi_{C}}C\in\Cyl(X)$.

  A \term{transverse graph} of a cylinder $C\in\Cyl(X)$ is an
  edge space or a vertex space of $C$. A transverse graph, when it
  doesn't matter if it's an edge space or vertex space, is typically
  denoted $F$. Choose an orientation on each edge of
  $\underlying(C)$ such that the edges of $\underlying(C)$ are
  $e_0,\dotsc,e_{n-1}$ and $\tau(e_i)=\iota(e_{i+1\bmod n}),$ and with
  vertices $v_i$ such that $\iota(e_i)=v_i$. Let $\alpha_{C}$ be
  the map 
  \begin{eqnarray*}
  \coprod_i\left(\tau_i\sqcup\iota_i^{-1}\right)\colon\edgegraphs(C)\sqcup\vertexgraphs(C)\to\\
  \hfill \vertexgraphs(C)\sqcup\edgegraphs(C)
  \end{eqnarray*}
  $\alpha_{C}$ respects the ordering and $\alpha_{C}^2$
  represents one $n$-th of a rotation of
  $C$. Also, $\alpha_{C}^{2n}=\mathrm{id}$.

%
  Let $X$ be a graph of spaces. If $X_1,\ldots, X_n$ are
  the irreducible components of $X,$ then each cylinder has image
  contained in one $X_i$. The \term{essential boundary},
  $\partial_{X}^{\mathrm{ess}}C$ is $\partial_{X_i}C$
  if $C$ has image in $X_i$.

  If $C$ is a cylinder of $X,$ $F$ a transverse graph of
  $C,$ and
  $\left|F\cap\partial_{X}^{\mathrm{ess}}(C)\right|>1,$ then
  the cylinder is \term{good}. Otherwise it is \term{bad}. Note that
  an irreducible component that consists of a single cylinder is
  automatically bad since $\compcomp$ of a cylinder is empty.
\end{definition}

\begin{definition}
  Let $\mathcal{F}(E)$ be the set of edge spaces
  \[\set{F\in\cup_{C\in\Cyl(X)}\edgegraphs(C)}:{\psi_{C}(F)\subset E}\]

  An element $F\in\mathcal{F}(E)$ is a \term{peripheral} element
  of $E$ if it contains a vertex $w\in\compcomp(X)\cap E$
  and if $w\in F'\in\mathcal{F}(E)$ then $F=F'$. The vertex $w$ is
  a \term{boundary vertex} of $E$.
\end{definition}

\par To get the ball rolling we need a way to take a peripheral
element $F$ of the set of cylinder cross sections $\mathcal{F}(E)$
and a boundary element $\compcomp(X)\cap F$ and push it around the
graph of spaces until a splitting vertex is discovered.

\begin{definition}[Pushing]
  A subset of a graph of spaces is \term{vertical} if it lies in a
  fiber of the map $\pi\colon X\to\underlying(X)$. Let
  $E\times\mathrm{I}$ be an edge space of $X$. Say that $x$
  and $y$ are equivalent if $x$ and $y$ have the same $E$
  coordinate. Horizontality is the equivalence relation generated by
  the relations on the edge spaces.

  If $Y$ is a vertical subset of $X$ then $Y$ pushes along a path
  $p\colon\left[0,a\right]\to\underlying(X)$ if there is a function
  $P$ such that for each $y\in Y,$ $P(y,\left[0,a\right])$ is
  horizontal, and the following diagram commutes
  \centerline{%
    \xymatrix{%
      Y\times\left[0,a\right]\ar[d]\ar[r]^{P(y,t)} &  X\ar[d]^{\pi} \\
      \left[0,a\right]\ar[r]^{p(t)} & \underlying(X)
    }
  }
  Since the diagram commutes, for each $t,$ the set $P(Y,t)$ is
  vertical.
\end{definition}

\par If $C$ is a cylinder of $X$ and $F$ is a transverse graph
of $C$ then $C$ is a vertical subset of $C$ and
$\psi_{C}(F)(=F)$ is a vertical subset of $X$. If a connected
vertical subset $Y$ of $X$ containing $F$ pushes along every
path $p$ that $F$ pushes along then $Y=F$. The rotation by $t$ of
$C$ is a one parameter family of homeomorphisms
$\alpha\colon C\times\mathbb{R}\to C$. If $F$ is a transverse
graph of $C$ then pushing $\psi_{C}(F)$ around $X$ can be
realized by the composition $P=\psi_{C}\circ\alpha$. 

\begin{lemma}[Separating Subgraphs, Structure of Vertex/Edge Spaces]
  \label{separatingsubgraphs}
  \par  Suppose $X$ is irreducible and separable.

  \par If $C$ is a cylinder of $X,$ and $F$ is a transverse
  subgraph of $C,$ then $\psi_{C}$ embeds $F$.  Every nonzero
  weight edge or vertex space of $X$ is a union of images of edge
  or vertex spaces, respectively, of cylinders of $X$.

  \par If $E$ is an edge space and $w\in\compcomp(X)\cap E$
  is not a cutpoint of $E$ then $\mathcal{F}(E)$ has a
  peripheral element containing $w$.
\end{lemma}

\begin{proof}
  \par Let $F$ be an edge space of $C$. Suppose there are vertices
  $p$ and $q$ such that $\psi_{C}(p)=\psi_{C}(q)$. First, note
  that $p$ and $q$ must be contained in
  $\partial_{X}(C)$. There are subgraphs $\Gamma_p$ and
  $\Gamma_q$ of $\Gamma(C)$ containing $p$ and $q,$
  respectively. Suppose $\psi_{C}(\alpha_{C}^k(p)) =
  \psi_{C}(\alpha_{C}^k(q))$ for all $k$. Then
  $\psi_{C}(\Gamma_p)$ and $\psi_{C}(\Gamma_q)$ must represent
  the same element of $\mathbb{S}(X),$ thus $\Gamma_p=\Gamma_q$ as
  sets, but this implies that $\Gamma_p$ must represent a periodic
  path in $\compcomp(X)$. The other possibility is that there
  exist $p$ and $q$ such that $\psi_{C}(p)\neq\psi_{C}(q),$
  but
  $\psi_{C}(\alpha_{C}(p))=\psi_{C}(\alpha_{C}(q))$. This
  is clearly impossible since edge maps of $X$ are injective. Thus
  $\psi_{C}$ embeds vertex and edge spaces. The collection
  $\coprod\psi_{C}$ is clearly injective on the collection of
  edges of vertex and edge spaces. Since every edge of an edge or
  vertex space comes from an annulus in $X,$ we have the first
  part of the lemma.
  
  \par Suppose $F,F'\in\mathcal{F}(E)$ with vertices $w_1,w_2\in
  F\cap F'$.  Clearly $w_1,w_2\in\compcomp(X)$. If $F$ is an edge
  graph of $C$ and $F'$ is an edge graph of $C',$ let
  $\Gamma_i^{(\prime)}$ be the component of $\compcomp(C^{(\prime)})$
  containing $\psi_{C^{(\prime)}}^{-1}(w_i)$. Since
  $\compcomp(X)\cap E$ doesn't separate $F$ and $F',$ we must
  have $\Gamma_i=\Gamma'_i,$ contradicting the construction of
  $\Cyl(X)$.
\end{proof}

\begin{figure}[tb]
  \psfrag{psi}{$\psi_{C}(E_{C})$}
  \psfrag{incompcomp}{$\in\compcomp(X)\cap E$}
  \centerline{%
    \includegraphics[scale=0.9]{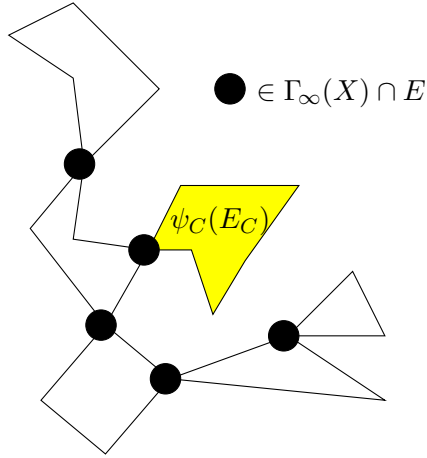}
  }
\caption{Edges are treelike. $\psi_{C}(E_{C})$ is the
image of an edge space $E_{C}$ of a cylinder $C$ of
$X$. The same picture holds for vertex spaces.}
\label{edgesaretreelike}
\end{figure}

We're now ready to prove Theorem~\ref{addroottofreegroup}. If we
represent maximal corank homomorphisms as immersions
$\Gamma(X)\to\underlying(X)$ factoring through
$\pi_1(X)\onto\pi_1(\underlying(X)),$ then we may choose an optimal
representation: by Lemma~\ref{foldtovalencethree} we may fold $X$ to
a space $X_c$ which minimizes $c$. Since $\chi(\underlying(X_c))
=\chi(\underlying(X)) =\chi(\Gamma(X)) =\chi(\Gamma(X_c)),$ by
Lemma~\ref{sameeulercharacteristic}, $X_c$ is separable. Also note
that the rest of the diagram in Figure~\ref{buildgraphofspaces}
commutes. $\psi$ is the homotopy equivalence given by the sequence of
folds and collapses to $X_c$.

\begin{figure}[ht]
\centerline{%
  \xymatrix{%
    \pi_1\Gamma_i\ar[r]^{\imath_*}\ar[d]^{\psi_*} & \pi_1 X_c\ar[r]\ar[d]^{\psi_*} & \pi_1\underlying(X_c)\ar@{>>}[d]^{{(\psi_{U})}_*} & \\
    \pi_1\Gamma_i\ar[r]^{\imath_*}   & \pi_1 X=G\ar[r]   & \pi_1\underlying(X)\ar@{>>}[dr]\ar@{>>}[r] & \pi_1 R_n\ar[d]^{\sim}\\
    F_i\ar[rrr]^{\phi}\ar[u] & & & \free'_n
  }
}
\caption{The commutative diagram representing $\phi$.}
\label{buildgraphofspaces}
\end{figure}

\begin{proof}[Proof of Theorem~\ref{addroottofreegroup}]
  Let $\phi:F=\free_n\into\free_n$ be a map that extends to a
  surjection $\widetilde{\phi} \colon  G=
  F\left[\set{\sqrt[k_i]{\gamma_i}}_{i=1..m}\right] \onto\free_n,$ with
  $\gamma_i$ pairwise nonconjugate, indivisible, and $k_i>1$. It's
  clear that $ G\onto\free_n$ is maximal corank. By the previous
  discussion, represent $\tilde{\phi}$ as a map of a separable space $
  X_c$ onto its underlying graph.

  Since $X_c$ is obtained from $X$ by a sequence of folds, every edge
  of $X$ is homeomorphic to some edge of $X_c,$ hence if the edges of
  $X_c$ are trees then so are those of $X$.

  \par Let $M_i$ be the mapping cylinder of the $k_i$--fold
  cover $S^1\to S^1$ corresponding to adding
  the $k_i$--th root to $\gamma_i$. The domain $S^1$ has an
  immersion $\gamma_i\colon S^1\to\compcomp(X_c)$ representing
  the conjugacy class of $\gamma_i$. Since $\gamma_i$ is indivisible,
  $\gamma_i$ is an element of $\mathbb{S}(X_c)$. The range $S^1$
  represents the $k_i$--th root of $\gamma_i$ and is called
  $r_i$. There is a map
  $\psi_{M_i}\colon M_i\to X_c$ which factors
  through some cylinder inclusion $\psi_{C}$. This map gives
  $M_i$ the structure of a 2-covered graph of spaces.
  
  \par First, note that since $\gamma_i$ is indivisible,
  $M_i$ embeds in the cylinder $C$. If $C$ was the
  union of more than one mapping cylinder, then some pair $\gamma_i$
  and $\gamma_j$ would have to be conjugate, thus $C=M_i$
  and our separable space $X_c$ is the union
  $\compcomp(X_c)\cup_{\gamma_i}M_i$.

  \par To complete the analysis of the cylinders, note that we must
  have $\vert F\cap\partial_{X_c}(M_i)\vert=k_i$ for any transverse
  graph $F,$ otherwise the immersion $\gamma_i$ must be a proper
  power. Thus every $F$ is a tree, and by
  Lemma~\ref{separatingsubgraphs}, every cycle in $E$ is contained in
  some element of $\mathcal{F}(E),$ $E$ is a tree.
\end{proof}

\begin{remark}
  We can now deduce Theorem~\ref{maintheorem} in the case that
  $\scott(G)=\scott(H)=(q-1,0)$. Since edge spaces of $X$ are trees,
  there is some element $\gamma_1$ which crosses an edge of
  $\Gamma(X)$ only once. Thus $G$ can be written as $F*\langle
  \gamma_1\rangle$.
\end{remark}

\section{Splitting graphs of spaces} 
\label{splittinggraphsofspaces}

\par Peripheral elements of $\mathcal{F}(E)$ and boundary vertices
play an essential role in finding moves which simplify graphs of
spaces. It is not enough that a graph of spaces merely have splittable
vertices. The notion that suffices since that of a splittable vertex
does not is that of a \term{splitting} vertex. A splitting vertex has
the property that one can collapse the ``outgoing'' edge adjacent to
the vertex, and strategically fold two of the edges in the resulting
graph of spaces, producing a new graph of spaces which \emph{still}
has a splitting vertex. Pushing and the treelike structure of edge
spaces are used to produce splitting vertices.

\begin{definition}[Splitting Vertex]
  A vertex $v$ of $\underlying(X)$ is \term{splitting}
  if 
  \begin{itemize}
  \item $v$ is splittable.
  \item $V\cong E_1\vee E_2$. The edge $e$ such that
    $E\cong V$ is the \term{outgoing} edge of $v$.
  \item If $w$ is the valence three vertex of $V$ and $e_1$ and
  $e_2$ are the other edges incident to $v,$ for at least one of
  $E_1$ or $E_2,$ say $E_1,$ there is a peripheral element
  $F\in\mathcal{F}(E_1)$ with boundary vertex $w\in F$ such that
  $\tau_1(w)$ is the valence three vertex of $\iota(e)$.
  \end{itemize}
  
  The edges $e_1$ and $e_2$ are the \term{incoming} edges. Numbered as
  in the bullets, $e_1$ is the \term{primary} incoming edge.
\end{definition}

\par The relationship between pushing, peripheral elements of
$\mathcal{F}(E),$ and splitting edges of $\underlying(X)$ is
what allows us to take a separable graph of spaces and convert it to
one with a bad cylinder. We first show that edge spaces have
peripheral elements.

\begin{lemma}
  \label{fixedpointfree}
  If $X$ is an irreducible, separable graph of spaces,
  $\chi(\underlying(X))<0,$ all of whose cylinders are good, then $X$
  has a splitting vertex.
\end{lemma}

\begin{proof}
  Let $\pi$ be the quotient map $X\to\underlying(X)$.

  Note that by Lemma~\ref{separatingsubgraphs} every edge space
  $E$ contains a vertex $w\in\compcomp(X)\cap E$ such that
  $w$ is contained in exactly one member $F$ of
  $\mathcal{F}(E)$. Choose such an edge $e$ not contained in
  $\cup\Gamma^i$ and regard $F,$ $w,$ as subsets of
  $E\times\set{1/2}$.

  Let $p\colon\left[0,a\right]\to\underlying(X)$ be the shortest
  path such that $p(0)=\pi(F)$ and $P(w,a)$ is a valence three vertex
  $w'$ of $\compcomp(X)$ in the vertex space $V$. Let
  $e_0,\dotsc,e_n$ be the sequence of edges that $p$ traverses. At
  integer values of $t,$ $F_t=P(F,t)\subset E_i\times\set{1/2}$. By
  the construction of the cylinders, $F_t\in\mathcal{F}(E_i),$ and
  it is obvious that $F_t$ is a peripheral element of the associated
  edge space and $P(w,t)$ is a boundary vertex of $E_i$ for the
  appropriate $t$. Since $P(w,a)$ is the valence three vertex of
  $V$. By construction $v$ is splitting.
\end{proof}

\par Now we need to know how to proceed when a separable graph of
graphs has a splitting vertex. There is a move, called \term{splitting},
which takes as input a separable graph of spaces which has a splitting
edge and outputs a ``simpler'' graph of spaces which is either
reducible or has a splitting edge and lower complexity.

\begin{definition}[Splitting]
  Suppose $X$ is separable and has a splitting vertex $V$. A
  splitting of $X$ is a graph of spaces $X_s$ obtained as follows:
  $V$ is splittable, so we can express $V$ as a wedge
  $V=L\vee_w R,$ with incident edge graphs homeomorphic to $L,$
  $R,$ or $V$.

  Define $e=e(v)$ to be the edge of $\underlying(X)$ such that
  $E\cong V$. Let $\overline{X}$ be the space obtained by collapsing
  $e$. Suppose $e_1$ and $e_2$ are the (oriented) edges other than $e$
  incident to $v$.  Let $v'$ be the other endpoint of $e$. Note that
  $v'$ is distinct from $v$ since
  $\mathscr{W}(V)>\mathscr{W}(L),\mathscr{W}(R)$.  Let $e_3$ and
  $e_4$ be the (oriented) edges other than $e$ incident to $v'$. In
  the collapsed space, let $\overline{e_i},$ $i=1,2,3,4,$ be the image
  of $e_i$. A splitting of $X$ is a nontrivial fold of
  $\overline{X}$ obtained by folding with $J=\set{1,3}$ or $J=\set{1,4}$.
\end{definition}

\par A splittable vertex $v$ of $X$ determines an edge
$e(v)\in\underlying(X)$ with $\iota(e)=v\neq\tau(e)=v'$. Let $w$
and $w'$ be the valence three vertices of $V$ and $V',$
respectively, and let $\pi\colon X\to\overline{X}$ be the map
which collapses the edge $e(v)$.

\par Let $v_i,$ $i=1,\dotsc,n,$ be the vertices of
$\underlying(X)$. The relative weight of a vertex $v_i$ is the
quantity \[\mathscr{W}(X\setminus v_i)=\sum_{j\neq
i}\mathscr{W}(V_j)\]

\begin{lemma}[Splitting Decreases Relative Weights]
  \label{relativeweights}
  \par If $v$ is a splitting vertex of an irreducible, separable,
  graph of spaces $X,$ then there is a collapse
  $\overline{X},$ a fold $X_s$ of $\overline{X}$ and if
  $X_s$ is irreducible, there is a splitting vertex $v_s$ of
  $X_s$ such that $\mathscr{W}(X\setminus
  v)>\mathscr{W}(X_s\setminus v_s)$.

  If $X$ is separable, then there is a sequence of collapses and
  folds to a space with no splitting vertices.
\end{lemma}

\begin{figure}[tb]
  \centerline{%
    \includegraphics{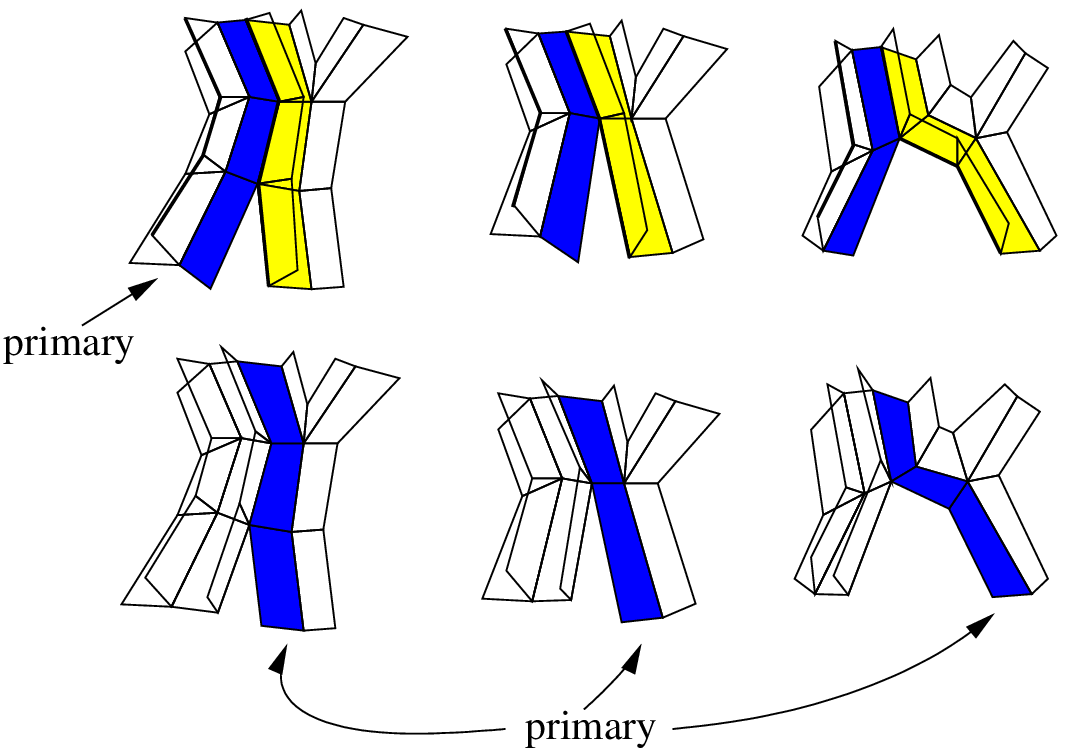}
  }
\caption{Illustration of case $\boldsymbol{\pi(w)\neq\pi(w')}$ of
Lemma~\ref{relativeweights}.}
\label{typesofsplitting}
\end{figure}

\begin{proof}
  Let $v=\tau(e),$ and let $g$ be the edge, not equal to $e,$ such
  that $V\cong E\vee G$. Let $v'$ be the terminal endpoint
  of $f,$ and let $h$ and $i$ be the two additional edges incident to
  $v'=\tau(f)$. Also, let $w$ be the separating vertex of $V$ and
  let $w'$ be the separating vertex of $V'$. Since
  $\mathscr{W}(F)>\mathscr{W}(E),\mathscr{W}(G),$ $f$ is
  embedded, thus we can collapse $f$ to obtain a space
  $\overline{X}$ with vertex $\overline{v},$
  $\overline{V}\cong V',$ and incident edge spaces $E,$
  $G,$ $H,$ and $I$.

  First, write $\overline{V}$ as
  $A\vee_{\overline{w}}B\vee_{\overline{w}}C$ such that
  $H\cong A\vee B$ and $I\cong B\vee C$. Let
  $\pi\colon X\to\overline{X}$ be the quotient map. There are
  two cases to consider.

  $\boldsymbol{\pi(w)=\pi(w')}:$ In this case, since $\pi(w)$
  separates, and $\pi(E)$ has only one element
  $F\in\mathcal{F}(E)$ such that $\pi(F)$ meets $\pi(w),$
  $\pi(E)$ is contained in, without loss, $A$. Folding
  $\overline{e}$ and $\overline{h}$ together creates two new vertices,
  one of which is homeomorphic to $H,$ is splittable, has an
  incident edge $e_s$ such that the pair $(e_s,h_s)$ is either
  splitting or such that $X_s$ has a weight $0$ edge. In the event
  that $(e_s,h_s)$ is splitting, the other vertex has weight
  $\mathscr{W}(V')-\mathscr{W}(E),$ i.e.,
  $\mathscr{W}(X_s\setminus\tau(e_s))<\mathscr{W}(X\setminus\tau(e))$.

  $\boldsymbol{\pi(w)\neq\pi(w')}:$ This case splits into two
  sub-cases. If $\pi(w')\subset\pi(G)$ then , without loss,
  $\pi(E)\subset\pi(H)$. Folding $\overline{h}$ and
  $\overline{e}$ together as in the previous case shows the lemma.

  We're left with the case $\pi(w')\subset\pi(E)$. Without loss,
  $\pi(G)\subset\pi(I)$. Folding $\overline{g}$ and
  $\overline{i}$ together creates a new splittable vertex with space
  isomorphic to $I$ and with incident edges isomorphic to $G$ and
  $B\cup(C\cap E)$. The vertex incident edge $g$ is splitting. This
  case is illustrated in the bottom row of
  Figure~\ref{typesofsplitting}.

  To see the second part of the lemma, suppose $X$ is irreducible
  and has a splitting pair. By the previous part of the lemma, we can
  split and fold to a space $X_s$ with a weight $0$ edge. Let
  $Y_i$ be the irreducible components of $X_s$. Each component
  $Y_i$ is seen to be separable. Now induct on
  $\chi(\underlying(Y_i))$. 
\end{proof}

\begin{figure}[tb]
  \centerline{%
    \includegraphics[scale=0.6]{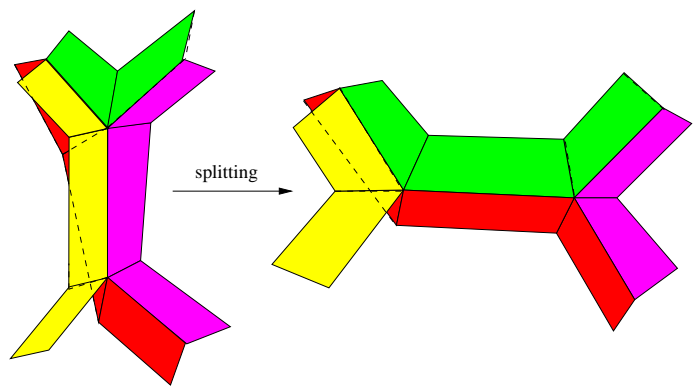}
    }
    \caption{Splitting when $\pi(w)=\pi(w')$.}
    \label{projectionsequal}
\end{figure}

\par If an irreducible component of a graph of spaces has a bad
cylinder then there is no guarantee the space can be further
simplified. The next theorem shows that one can convert a graph of
spaces to a ``minimal'' one, where minimal means that no sequence of
collapses and folds ever leads to the creation of a weight $0$
edge. Combining Lemmas~\ref{fixedpointfree} and~\ref{relativeweights}
we have the following theorem.

\begin{theorem}[Splitting to bad cylinders]
  \label{splittocylinders}
  If $X$ is separable, there is a sequence of splittings to a
  space $X_b$ such that every irreducible component has a bad
  cylinder.
\end{theorem}

As a consequence of this and the analysis of edge spaces from the
previous section, we can now establish the conjugacy separability
result stated in the introduction.

\begin{figure}[tb]
  \psfrag{gamma3}{$\gamma^1_j\in\mathrm{Z}_i$}
  \psfrag{gamma4}{$\gamma^2_j\in\mathrm{Z}_{i'}$}
  \psfrag{gamma1}{}
  \psfrag{gamma2}{}
  \psfrag{MI}{$A_j\mbox{from }t_j$}
  \psfrag{C}{$C$}
  \psfrag{F}{$E$}
  \centerline{%
    \includegraphics[scale=0.6]{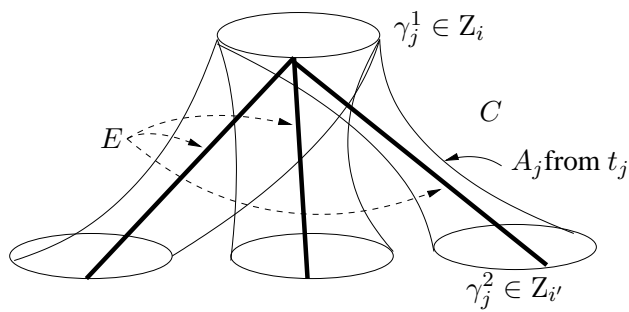}
  }
  \caption{A typical cylinder from Corollary~\ref{maximumrankconjugacy}.}
  \label{hnnfigure}
\end{figure}

\begin{proof}[Proof of Corollary~\ref{maximumrankconjugacy}]
  As before Theorem~\ref{addroottofreegroup}, represent the
  homomorphism $\widetilde{\phi}\colon G\onto\free_n$ as a
  homomorphism $X_c\to\underlying(X_c)$. Since $\widetilde{\phi}$ has
  maximal corank, $X_c$ is separable.

  We prove the theorem by observing that the hypothesis that $\sim$
  has no singleton equivalence classes implies that either all
  cylinders are good or the theorem holds. What are the cylinders
  of $X_c$? The maximal abelian subgroups $\mathrm{Z}_i$ of $F$
  can be represented as elements of $\mathbb{S}(X_c)$. The stable
  letters $t_j$ from $G$ give $\gamma^1_j\in\mathrm{Z}_i$ and
  $\gamma^2_j\in\mathrm{Z}_{i'},$ and, for each $j,$ an annulus
  $A_j$ glued between $\mathrm{Z}_i$ and $\mathrm{Z}_{i'}$ as
  elements of $\mathbb{S}(X_c)$. Then the cylinders of $X_c$
  are represented precisely by the equivalence classes of $\sim$ from
  the statement of the theorem. Since there are no singleton
  equivalence classes the boundary $\partial_{X_c}(C)$ of
  every cylinder $C$ has more than one component. The key thing to
  notice is that an edge space $E$ of $C$ meets \emph{every}
  component of $\partial_{X_c}(C)$ \emph{at least} one
  time. Thus a cylinder is bad if and only if
  $\partial_{X_c}^{ess}(C)$ has only one component. A cylinder
  is illustrated in Figure~\ref{hnnfigure}.

  By Theorem~\ref{splittocylinders} we may replace $X_s$ by a
  separable graph of spaces $X_b$ whose irreducible components
  each contain a bad cylinder.  Choose a bad cylinder $C$ and a
  component
  $\mathrm{Z}_b\subset\partial(C)\setminus\partial^{ess}(C)$. All
  edges of $\compcomp(X_b)$ which meet $\mathrm{Z}_b$ have weight
  $0$. This collection of edges can be folded together to give $F$ a
  free factorization $\mathrm{Z}_b*F'$ satisfying the theorem.
\end{proof}

\bibliographystyle{amsalpha}
\bibliography{krull}

\end{document}